\documentclass[11pt]{article}

\usepackage{arxiv}
\usepackage{amsmath}
\usepackage{graphicx}
\usepackage{enumerate}
\usepackage{natbib} 
\usepackage{hyperref}
\usepackage{xcolor}
\hypersetup{
    hidelinks
}

\usepackage[T1]{fontenc}
\usepackage{anyfontsize}

\usepackage{glossaries}
\usepackage{placeins}
\newcommand{\ra}[1]{\renewcommand{\arraystretch}{#1}}
\usepackage{booktabs}
\newcommand{\A}{\mathcal{A}}
\newcommand{\B}{\mathcal{B}}

\newcommand{\D}{\mathcal{D}}

\newcommand{\X}{\mathcal{X}}

\newcommand{\V}{\mathcal{V}}

\newcommand{\exclude}[1]{}{}
\newcommand{\taur}{t^\text{res}} 
\newcommand{\tauw}{t^\text{wai}} 
\newcommand{\taud}{t^\text{det}} 
\newcommand{\taup}{t^{p}} 
\newcommand{\taua}{t_{a}} 
\newcommand{\decs}{d^{\textnormal{sngl}}}
\newcommand{\decm}{d^{\textnormal{mlti}}}
\newcommand{\decp}{d^{\textnormal{pool}}}
\newcommand{\decq}{d^{\textnormal{queu}}}
\newcommand{\decrel}{d^{\textnormal{relc}}}
\newcommand{\decrec}{d^{\textnormal{rchr}}}
\newcommand{\decidle}{d^{\textnormal{idle}}}
\newcommand{\deccont}{d^{\textnormal{cont}}}
\newcommand{\Decs}{\mathcal{D}^{\textnormal{sngl}}}
\newcommand{\Decm}{\mathcal{D}^{\textnormal{mlti}}}
\newcommand{\Decp}{\mathcal{D}^{\textnormal{pool}}}
\newcommand{\Decq}{\mathcal{D}^{\textnormal{queu}}}
\newcommand{\Decrel}{\mathcal{D}^{\textnormal{relc}}}
\newcommand{\Decrec}{\mathcal{D}^{\textnormal{rchr}}}
\newcommand{\Decidle}{\mathcal{D}^{\textnormal{idle}}}
\newcommand{\Deccont}{\mathcal{D}^{\textnormal{cont}}}

\newcommand{\Aempt}{\A^{\textnormal{empt}}}
\newcommand{\Aoccu}{\A^{\textnormal{occu}}}

\newcommand{\detourPenalty}{\rho^\text{detr}}
\newcommand{\rechargeCost}{\rho^\text{rchr}}
\newcommand{\myopicpolicy}{\pi^\text{myo}}

\newcommand{\pathend}{p_\text{end}}
\usepackage{mathtools}

\usepackage{amsfonts,amssymb,array}
\usepackage{amsthm}

\usepackage{algorithm}
\usepackage{algorithmic}

\usepackage{cleveref}
\usepackage{siunitx}
\usepackage{multirow}
\usepackage{multicol}
\usepackage{subcaption}
\usepackage{tikz}
\usepackage{fontawesome5}
\usepackage{bm}
\usepackage{enumitem}

\DeclareMathOperator*{\argmax}{arg\,max}

\newtheorem{dfn}{Definition}[section]
\newtheorem{proposition}{Proposition}[section]

\newcommand{\rangeICE}{\qty{26}{\hour}}
\newcommand{\rangeEV}{\qty{17}{\hour}\,\qty{41}{\minute}}
\newcommand{\rechargingRateICE}{\qty[per-mode=symbol]{2.308}{\second\per\hour}}
\newcommand{\rechargingRateDCFC}{\qty[per-mode=symbol]{2.262}{\minute\per\hour}}
\newcommand{\rechargingRateLTwoC}{\qty[per-mode=symbol]{45.23}{\minute\per\hour}}
\newcommand{\symbolRechargingRate}{\dot{\ell}} 

\newcommand{\emptycar}{\resizebox{!}{0.85cm}{\color{black}\faTaxi}}
\newcommand{\occupiedcar}{\resizebox{!}{0.85cm}{\color{blue!50}\faTaxi}}

\newcommand{\originnodem}{\resizebox{!}{0.8cm}{\color{black}\faMale}}

\newcommand{\destinationnode}{\resizebox{!}{0.6cm}{\color{black}\faHome}}
\newcommand{\currentlocation}{\resizebox{!}{0.6cm}{\color{black}\faMapMarker*}}
\newcommand{\pastlocation}{\resizebox{!}{0.6cm}{\color{gray}\faMapMarker*}}

\newcommand{\tInSetTimeSteps}{t\in\{1,\dots, T\}}

\newcommand{\partialOrderDefinition}[1]{\begin{dfn}[Partial Order]\label{def#1:partial_order}
        Let $a, a'\in\A$ with $a = (o_{a} \; d_{a} \; l_{a} \; n_{a} \; t_{a})^{\intercal}$ and $a' = (o_{a'} \; d_{a'} \; l_{a'} \; n_{a'} \; t_{a'})^{\intercal}$. Then, $a \preccurlyeq a'$ if, and only if, $(o_{a},d_{a}) = (o_{a'}, d_{a'})$, $l_{a} \leq l_{a'}$, $n_{a} \leq n_{a'}$ and $t_{a} \geq t_{a'}$.
    \end{dfn}
}

\newcommand{\monotonicityProposition}[1]{
\begin{proposition}[Monotonicity]\label{prop#1:monotonicity}
    For every pair of vehicle attributes $a,a' \in \A$ with $a\preccurlyeq a'$, if the statements
    $a^M(a,d) \preccurlyeq a^M(a',d)$ and $c_{tad} \leq c_{ta'd}$
    are true for all decisions $d\in \D(a) \cap \D(a')$ and time steps  $\tInSetTimeSteps$, then the following monotonic properties hold:
    \begin{enumerate}
        \item ${v}^{R*}_a$ increases monotonically with the $l_a$ dimension: $l_{a} < l_{a'} \implies {v}^{R*}_a \leq {v}^{R*}_{a'}$.
        \item ${v}^{R*}_a$ increases monotonically with the $n_a$ dimension: $n_{a} < n_{a'} \implies {v}^{R*}_a \leq {v}^{R*}_{a'}$.
        \item ${v}^{R*}_a$ decreases monotonically with the $t_a$ dimension: $t_{a} > t_{a'} \implies {v}^{R*}_a \leq {v}^{R*}_{a'}$.
    \end{enumerate}
\end{proposition}
}

\newcommand{\detourPenaltyMonoDefinition}[1]{\begin{dfn}[Monotonicity of detour penalty]
        \label{def#1:monotonicity_of_detour_penalty}
        The detour penalty function $\detourPenalty{}$ is \emph{monotonically increasing in the path length} if, and only if, for all vehicle attributes $a\in\Aempt$, all sets of requests $B\subseteq\mathcal{B}$ with $|P_{aB}| > 0$, and all pairs of paths $p, p' \in P_{aB}$, the following statement is true:\[
        \tau(p) \leq \tau(p') \implies \detourPenalty{}(a,\decm_{Bp}) \leq \detourPenalty{}(a,\decm_{Bp'}).
        \]
\end{dfn}}

\newcommand{\enumerationDefinitions}[1]{
    \begin{dfn}[Sufficient paths set]\label{def#1:sufficientset}
        Let $\pi$ be a policy. For all time steps $\tInSetTimeSteps$ and all states $S_t$ where a vehicle attribute $a\in\Aempt$ satisfies $R_{ta} > 0$ and where the nonempty set of requests $B\subseteq\mathcal{B}$ satisfies $\forall b\in B.D_{tb} > 0$ and $|P_{aB}| > 0$, the set of paths $P^*_{aB}\subseteq P_{aB}$ is \emph{a-B-$\pi$-sufficient} if, and only if, there exists a solution $x_t$ in the set of optimal solutions $\pi(S_t)$ such that all decision variables $x_{tad}$ with $d\in \{\decm_{Bp} \in \Decm_a \mid p\in P_{aB}\setminus P^*_{aB}\}$ are equal to zero.
    \end{dfn}
    \begin{dfn}[\Gls{ldsps}]\label{def#1:shortestpermutationpathsset}
        For all time steps $\tInSetTimeSteps$ and all states $S_t$ where a vehicle attribute $a\in\Aempt$ satisfies $R_{ta} > 0$ and where the nonempty set of requests $B\subseteq\mathcal{B}$ satisfies $\forall b\in B.D_{tb} > 0$ and $|P_{aB}| > 0$, the set of paths $P_{aB}^* \subseteq P_{aB}$ is an \emph{\gls{ldsps} for (a,B)} if, and only if the following statement is true for all requests $b\in B$: Either $P_{aB}^*$ contains exactly one of the shortest paths in $P_{aB}$ that end in the destination of $b$, $d_b$, or no path in $P_{aB}$ ends in $d_b$.
\end{dfn}}

\newcommand{\enumerationProposition}[1]{\begin{proposition}[Sufficiency of the \gls{ldsps}]\label{prop#1:sufficientset}
        For all time steps $\tInSetTimeSteps$ and all states $S_t$ where a vehicle attribute $a\in\Aempt$ satisfies $R_{ta} > 0$ and where the nonempty set of requests $B\subseteq\mathcal{B}$ satisfies $\forall b\in B.D_{tb} > 0$ and $|P_{aB}| > 0$, if the approximate value function used by policy $\pi^\mathrm{VFA}$ satisfies \Cref{prop:monotonicity} and $\detourPenalty{}$ is monotonically increasing in the path length, all \glspl{ldsps} for $(a,B)$ are a-B-$\pi^\mathrm{VFA}$-sufficient.
\end{proposition}}
          \newcommand{\VFARfrMedianPoolingtrue}{89.2\%} \newcommand{\VFARfrMedianPoolingfalse}{86.4\%}                    \newcommand{\VFARfrIqrPoolingtrue}{3.9\%} \newcommand{\VFARfrIqrPoolingfalse}{6.9\%}                                                                                  \newcommand{\PMRfrMedianIqr}{76.8\% (24.2\%)}                                    \newcommand{\VFARfrMedianIqr}{88.4\% (4.9\%)} \newcommand{\VFARfrMedianIqrDCFC}{89.6\% (2.3\%)} \newcommand{\VFARfrMedianIqrLTwoC}{84.3\% (5.7\%)} \newcommand{\VFARfrMedianIqrICE}{89.2\% (3.7\%)}           \newcommand{\VFARfrMedianIqrPoolingtrueDCFC}{90.4\% (2.1\%)} \newcommand{\VFARfrMedianIqrPoolingfalseDCFC}{89.1\% (2.4\%)}     \newcommand{\VFARfrMedianIqrPoolingtrueLTwoC}{86.4\% (2.3\%)} \newcommand{\VFARfrMedianIqrPoolingfalseLTwoC}{81.3\% (5.0\%)}

\newcommand*\nestedglsentry[1]{%
  \protect\ifglsused{#1}{%
    \glsentryshort{#1}%
  }{%
    \glsentrylong{#1}%
  }%
}

\newacronym{ann}{ANN}{artificial neural network}
\newacronym{adp}{ADP}{Approximate Dynamic Programming}
\newacronym{darp}{DARP}{dial-a-ride problem}
\newacronym{mdp}{MDP}{Markov decision process}
\newacronym{lp}{LP}{Linear Program}
\newacronym{ip}{IP}{Integer Program}
\newacronym{pm}{PM}{parameterized myopic}
\newacronym{sdp}{SDP}{sequential decision-making process}
\newacronym{sdvrp}{SDVRP}{Stochastic Dynamic Vehicle Routing Problem}
\newacronym{rhs}{RHS}{ride-hailing system}
\newacronym{rrhs}{RRHS}{robotaxi ride-hailing system}
\newacronym{prrhs}{PRRHS}{pooling-enabled robotaxi ride-hailing system}
\newacronym{vfa}{VFA}{value function approximation}
\newacronym{NYC}{NYC}{New York City}
\newacronym{TLC}{TLC}{\nestedglsentry{NYC} Taxi and Limousine Commision}
\newacronym{ICE}{ICE}{internal combustion engine}
\newacronym{EV}{EV}{electric vehicle}
\newacronym{DCFC}{DCFC}{direct current fast charging}
\newacronym{l2c}{L2C}{level 2 charging}
\newacronym{rfr}{RFR}{reward fulfillment ratio}
\newacronym{iqr}{IQR}{interquartile range}
\newacronym{me}{MOE}{margin of error}
\newacronym{ldsps}{LDSPS}{last dropoff shortest paths set}
\newacronym{spprc}{SPPRC}{shortest path problem with resource constraints}

\title{Sequential Decision-making for Ride-hailing Fleet Control: A Unifying Perspective}

\author{
Stefan Pilot \\
  RWTH Aachen University\\
  \href{mailto:pilot@cl.rwth-aachen.de}{\texttt{pilot@cl.rwth-aachen.de}}
   \And
 Murwan Siddig \\
  University of Florida\\
  \href{mailto:msiddig@ufl.edu}{\texttt{msiddig@ufl.edu}}
}
\begin{document}
\maketitle

\begin{abstract}
  This paper provides a unified framework for the problem of controlling a fleet of ride-hailing vehicles under stochastic demand. 
  We introduce a sequential decision-making model that consolidates several problem characteristics and can be easily extended to include additional characteristics. To solve the problem, we design an efficient procedure for enumerating all feasible vehicle-to-request assignments, and we introduce scalable techniques to deal with the exploration-exploitation tradeoff. We construct reusable benchmark instances that are based on real-world data and that capture a range of spatial structures and demand distributions. Our proposed modelling framework, policies and benchmark instances allow us to analyze interactions between problem characteristics that were not previously studied. We find no significant difference between revenue generated by internal combustion engine fleets and fast-charging electric fleets, but both significantly outperform slow-charging electric fleets. We also find that pooling increases the revenue, and reduces revenue variability, for all fleet types. Our contributions can help coordinate the significant research effort that this problem continues to receive.
\end{abstract}


\section{Introduction}\label{sec:introduction}

Ride-hailing is an on-demand transportation service in which customers request rides, typically through a mobile application or an online platform.
The service provider matches the requests with available vehicles that transport the customers to their destinations. 
Ride-hailing service providers can also offer pooling options, where multiple customers share a ride at the same time, usually at a lower cost. 
Examples of ride-hailing services include Uber, Lyft, DiDi, and Waymo. 
Ride-hailing services have become increasingly popular in recent years due to their convenience, affordability, and flexibility compared to traditional taxi services and public transportation. They also have the potential to reduce traffic congestion and emissions by optimizing vehicle usage and reducing the number of empty seats on the road \citep{ackermann2023ejtl}.

Optimal control of \glspl{rhs} is a challenging sequential decision-making problem. Every decision made by the service provider impacts not only the current \gls{rhs} performance, but also the \gls{rhs}'s future performance. There are many possible decisions that the provider can make, e.g., dispatching vehicles to satisfy customer requests, keeping vehicles parked at their current locations, repositioning empty vehicles to different locations, and recharging vehicles if they are electric. The decision of dispatching a vehicle to satisfy customer requests is, in itself, a complex decision with many layers. For example, which vehicle should be assigned to which requests? If a vehicle is assigned to multiple requests, should the requests be pooled together or served one after the other? In what order should they be picked up and dropped off? All of these decisions are interrelated and must be executed dynamically while satisfying operational constraints (e.g., vehicle availability, capacity, driving range) and customer satisfaction constraints (e.g., waiting time until the vehicle arrives, detour time to pick up other passengers).

The complex nature of this optimal control problem has led to a rich and diverse body of literature studying various problem settings, modeling choices, and solution methods~\citep[see surveys][]{mourad:2019, soeffker:2022, HILDEBRANDT2023106071}. As a result, the literature appears fragmented and would benefit significantly from a more coordinated research effort that facilitates meaningful comparisons and generalization of findings across studies. This work is an attempt towards a unified framework for studying \gls{rhs} control.


\subsection{Scope and Terminology}\label{sec:terminology_and_scope}
\glspl{rhs} can take different forms depending on the service provider's business model, the type of vehicles used, and the level of automation. For instance, the fleet of vehicles used by a \gls{rhs} can consist of \gls{ICE} vehicles, \glspl{EV}, or a mix of both. We consider both \gls{ICE} vehicles and \glspl{EV}. Furthermore, \glspl{rhs} can be fully autonomous (e.g., Waymo and Zoox) or operated by human drivers. In the latter case, drivers can either be employees of the service provider (e.g., Ioki, Moia, and Via) or independent contractors (e.g., Uber and Lyft). We focus on \glspl{rhs} that operate a fleet of vehicles owned by the service provider, regardless of whether the vehicles are autonomous or driven by employees.

Different studies use different terminology to refer to the \glspl{rhs} of our interest \citep[see e.g.,][]{qin:2022}. The terminology is especially inconsistent in studies that allow the option of pooling.
We emphasize that the focus of this work does not include rental car sharing \citep{boyaci:2015}, sharing rides with parcels \citep{li:2014}, matching private drivers and riders with similar plans \citep{agatz:2012}\exclude{, auction-based ride assignments \citep{kleiner:2011}}, and dial-a-ride paratransit with large-capacity vehicles and pre-specified requests \citep{heitmann2024accelerating}. To that end, we use the term \emph{pooling-enabled} \gls{rhs} to emphasize the taxi-like on-demand nature of \glspl{rhs} with the option to share rides that we study.


\subsection{Problem statement}\label{sec:problem_statement}
We study a pooling-enabled \gls{rhs} that controls a fixed-size fleet of homogeneous vehicles. The fleet consists of either \gls{ICE} vehicles or \glspl{EV}, and each vehicle has a fixed capacity for passengers and a fixed maximum driving range before it needs to recharge or refuel. 

The \gls{rhs} dispatches vehicles to serve travel requests. 
Travel requests are not known in advance and arrive over the planning horizon. The specific point in time at which the operator makes decisions is called the decision epoch. At each decision epoch, the \gls{rhs} chooses which requests to accept, which vehicles to dispatch to serve the accepted requests, and how to dispatch them. Vehicles not serving requests may be relocated, left idle, or refueled (recharged, if \glspl{EV}).

The suitability of a vehicle to serve a request depends on the request characteristics and the vehicle state. A request is characterized by its origin, destination, headcount, latest response time, latest pickup time, and fare. If a request is not accepted in one epoch, it is reconsidered in the next until its latest response time passes, after which it is treated as rejected. The \gls{rhs} may pool accepted requests and serve them simultaneously with one vehicle. Pooling can increase travel time relative to serving requests individually, but passengers pay a lower fare. Pooling can also reduce waiting times in the rush hour if it allows the same fleet to serve more requests in a given period. The \gls{rhs} decides whether to pool and which vehicles to dispatch.

A vehicle state includes its current location, destination (if occupied), remaining capacity (empty seats), driving range (for \glspl{EV}, the battery level), and actionable time (when it can start a new task). A vehicle can be assigned to a request only if it has enough empty seats to accommodate the request headcount, can reach the request origin before the latest pickup time, and has sufficient driving range to reach the origin and then the destination. 

The \gls{rhs} aims to maximize its total expected reward over the planning horizon. The total reward is the sum of the fares collected from the passengers minus detour penalties and refuelling or recharging costs incurred by the system.


\subsection{Contributions}\label{sec:contributions}
The main contributions of this paper are:
\begin{itemize}
    \item \textbf{Model:} We present a unified sequential decision-making model that (1) consolidates many of the previously studied features, decisions, and constraints, and (2) can be easily extended to include additional features, decisions, constraints, and objectives.
    \item \textbf{Policies:} We present a tractable and well-performing \gls{vfa} policy that combines the advantages of existing approaches. This involves introducing an efficient procedure for enumerating all feasible vehicle-to-request assignments and scalable techniques for handling the exploration-exploitation tradeoff.
    \item \textbf{Benchmark instances:} We create a set of benchmark instances that (1) are based on real-world data, (2) include a variety of spatial structures and demand patterns, and (3) are easily accessible, and can be directly used for future research.
    \item \textbf{Numerical results:} We perform extensive numerical experiments, which show that our \gls{vfa} policy significantly outperforms a baseline policy, which is consistent with previous literature. We also analyze previously unstudied interactions between problem characteristics, and find that both \gls{ICE} fleets and fast-charging \gls{EV} fleets obtain high rewards compared to slow-charging \gls{EV} fleets.
    We also find that pooling increases the revenue of the \gls{rhs} and reduces revenue variability for all fleet types.
\end{itemize}

The remainder of this paper is organized as follows.
Section \ref{sec:related_work} surveys relevant literature.
Section \ref{sec:formulation} describes the sequential decision-making model.
Section \ref{sec:policies} defines our two policies.
Section \ref{sec:experiments} presents our benchmark instances and numerical experiments.
Section \ref{sec:conclusion} concludes the paper and discusses future research directions.
\FloatBarrier
\section{Related Work}\label{sec:related_work}

In this section, we provide an overview of the literature on the topic by highlighting the different decisions, constraints, and performance metrics that have been considered. We also discuss the different models and solution methods that have been used to solve the problem. 

\subsection{Decisions}\label{subsec:decisions2}
The most basic decision that the operator can make is to assign an \textit{empty} vehicle to a single travel request. We refer to this decision as the \textit{single trip} decision ($\decs$). The operator can also assign an \textit{empty} vehicle to park at its current location \citep{alkanj:2020,kullman:2022,yu2023coordinating}. We refer to this as the \textit{idle} decision ($\decidle$). Not all studies consider the $\decidle$ decision. For instance, \cite{alonso:2017} assume that empty vehicles are relocated to areas with high demand. We use a separate decision for modeling the relocation of empty vehicles. We refer to this decision as the \textit{relocate} decision ($\decrel$). The $\decrel$ decision is commonly studied in the literature \citep[e.g.,][]{alkanj:2020,kullman:2022,tuncel2023integrated}.

Some decisions are only relevant for certain types of problem settings. For instance, the \textit{recharge} decision ($\decrec$) is typically only included in studies that consider electric vehicles  \citep[e.g.,][]{alkanj:2020,kullman:2022,yu2023coordinating}. Similarly, two decisions that are only relevant for studies that allow ride pooling are the \textit{multi-trip} decision ($\decm$) and the \textit{pool} decision ($\decp$). The $\decm$ decision assigns an \textit{empty} vehicle to multiple requests, where the passengers of those requests may or may not occupy the vehicle simultaneously at some point during the trip \citep{alonso:2017,tuncel2023integrated}. The $\decp$ decision assigns an \textit{occupied} vehicle to a new request, which is picked up before the vehicle drops off its existing occupants at their destination \citep{yu:2020,heitmann2023combining}. To distinguish between $\decp$ and the decision of assigning an \textit{occupied} vehicle to a new request whose passengers are picked up \textit{after} the vehicle drops off its existing occupants, we introduce a separate decision that we refer to as the \textit{queue} decision ($\decq$). The $\decq$ decision is not always considered in the literature and is sometimes modeled implicitly as a $\decs$ decision \citep{kullman:2022}, or using tentative routes that can be adapted over time \citep{heitmann2023combining}.

Some decisions are only relevant for certain types of models. For instance, the specific point in time at which the operator makes decisions is called the decision epoch. Decision epochs can be triggered by events, e.g., arrival of a new request \citep{kullman:2022}, or time-triggered at regular intervals \citep{alonso:2017,alkanj:2020}. We consider a time-triggered decision epoch, but our model can be easily adapted to consider an event-triggered decision epoch. Modeling time-triggered decision epochs requires special handling of vehicles that are still performing tasks at the time of the decision epoch.
We do this by introducing a pseudo decision that we refer to as the \textit{continue} decision ($\deccont$), which ensures that the vehicle continues to perform the task it is currently performing (e.g., servicing requests, relocating or recharging). The $\deccont$ decision is referred as the \textit{null} decision by \cite{kullman:2022}.
 
\subsection{Constraints}\label{subsec:constraints}
The most commonly studied type of constraints concern the \textit{waiting} time ($\tauw$), which refers to the time elapsed between the moment at which a customer is assigned a vehicle and the moment at which the vehicle picks up the customer; the \textit{detour} time ($\taud$), which refers to the additional time incurred during one passenger's trip due to picking up or dropping off other passengers along the way; and the \textit{response} time ($\taur$), which refers to the time elapsed between the moment at which a customer sends a request for a ride and when the system assigns a vehicle to the customer. 
For instance, \citet{ozkan:2020} and \citet{kullman:2022} impose a maximum waiting time constraint to ensure that customers are picked up within a certain time after being assigned a vehicle. Some studies do not impose a constraint on the waiting time but penalize it in the objective function \citep{yu:2020}.
A detour time constraint is only relevant for studies that consider pooling-enabled \glspl{rhs}. For instance, \citet{alonso:2017} consider a detour time constraint implicitly by imposing a maximum drop-off time for every request and penalizing the detour time in the objective function. \cite{yu:2020} do not impose a detour time constraint but penalize it in the objective function.
Constraints on the response time \citep{erdmann2021combining} are not as common as the waiting and detour time constraints. However, a key benefit of modeling a response time constraint is that it defines a demand backlog for the requests that the system receives. Studies without a response time constraint typically consider all requests that are not served instantly to be lost \citep[e.g.,][]{alkanj:2020} or penalize the lost requests in the objective function \citep{yu:2020,alonso:2017}.
Another type of constraint concerns the recharge and idle decisions. While some studies assume that vehicles can charge and idle everywhere \citep{alkanj:2020}, others assume that vehicles can only charge in specific locations \citep{yu2023coordinating}, only idle in specific parking lots \citep{ackermann2023ejtl}, or both \citep{kullman:2022}. 

In this paper, we impose the waiting time and response time as hard constraints and penalize the detour time in the objective function. Our reasoning for this is as follows. Customers who share rides with other passengers typically pay cheaper fares, and it is reasonable to interpret a penalty cost in the objective function as a discount or compensation for the inconvenience of a detour time. Additionally, imposing the waiting time as a hard constraint allows us to restrict the time it takes a vehicle to pickup a request to some upper bound, e.g., the time it takes to complete the trip itself. Similarly, imposing the response time as a hard constraint allows us to model a demand backlog for the system that can be interpreted as what \citet{yan:2020} refer to as a ``batching window''. Furthermore, we assume that vehicles can idle and fully charge everywhere.
We make this assumption to simplify the problem. Nevertheless, our model can be easily adapted to consider dedicated charging stations and parking lots, as well as any combination of hard and soft constraints on the waiting time, response time, and detour time.

\subsection{Performance metrics}\label{subsec:performance_metrics}
The choice between imposing a constraint or adding a penalty to the objective function is subjective, in part because customers' tolerance for waiting, response, and detour times can vary significantly. Service providers try to strike a balance between profitability and customer satisfaction. 
Many performance metrics have been considered in the literature to calibrate this trade-off. The most commonly studied performance metrics are the rewards generated by the system for serving requests \citep[][]{alkanj:2020}, operational costs such as the cost of driving \citep{kullman:2022}, charging \citep{alkanj:2020}, as well as service quality metrics such as the waiting time \citep{alonso:2017,lyu2024multiobjective}, response time \citep{erdmann2021combining}, and detour time \citep{alonso:2017}. Some studies consider the number of served requests \citep{ozkan:2020}, cost of losing requests \citep{alonso:2017,yu:2020}, the number of available vehicles \citep{braverman:2019}, or geographical unfairness where customers whose origins and/or destinations are far from centralized locations suffer from excessive service rejections \citep{ben2025data}. Many of these performance metrics are often considered simultaneously. This could be as simple as maximizing a revenue function that comprises rewards and costs \citep{alkanj:2020,kullman:2022}, but can also be modelled as a multi-objective problem \citep{lyu2024multiobjective}. In this paper, we use a simple revenue function that comprises the trip fares (gained by serving requests) minus the detour penalty costs and the charging costs. This revenue function can be easily adapted to consider other performance metrics. Service fairness metrics are outside the scope of this paper.

\subsection{Models and solution methods}\label{subsec:models_methods}
The solution methods used in the literature can be distinguished based on the type of policy they study, and whether constructing the policy involves solving an assignment problem.
Most policies solve some variant of a linear assignment problem to match vehicles to requests or to any other task (e.g., $\decidle$ and $\decrel$). The assignment problem can be solved as a \gls{lp} \citep[e.g.,][]{alkanj:2020} or an \gls{ip} \citep[e.g.,][]{alonso:2017}.
There are also policies that do not solve an assignment problem. For example, \cite{kullman:2022} learn Q-value approximations using an \gls{ann} to make separate, decentralized decisions for each vehicle. 

Some studies use myopic policies. For instance, the seminal paper by \cite{alonso:2017} define a bipartite request-trip-vehicle (RTV) graph whose nodes represent vehicles and potential requests groups (i.e., shared rides) and whose edges represent feasible matches between vehicles and these rides. The authors then solve an assignment problem myopically, in a rolling-horizon fashion, without explicitly considering the downstream impact of the current decisions on the future performance of the system. A myopic policy can also be parameterized \citep{powell:2011}. For instance, \cite{alkanj:2020} solve an assignment problem that includes a parameterized rule to ensure that a vehicle is assigned to charge if its battery level is below a certain threshold.
Other studies create policies that explicitly consider the downstream impact of the current decisions by optimizing the current reward (or cost) plus an approximation of the value of future states, which the system evolves to, given the current decisions. We refer to these policies as \gls{vfa}-based policies. \gls{vfa}-based policies can be further distinguished based on the method used to approximate the value of future states. For instance, \cite{alkanj:2020} use lookup tables to approximate the value of future states
and \cite{wang:2018} learn an \gls{ann}-\gls{vfa} from the perspective of a single vehicle.
We refer the reader to \cite{qin:2022} for a more detailed overview of solution methods used in the literature.

In this paper, we study a \gls{pm} policy and a \gls{vfa} policy.
To design each policy, we use a modified version of the RTV graph algorithm to enumerate the decisions of a linear assignment problem. For the \gls{pm} policy, we solve the assignment problem as an \gls{ip} to obtain integer solutions directly. For the \gls{vfa} policy, we first solve a relaxed version of the assignment problem as an \gls{lp} to obtain a \gls{vfa} lookup table, which we then use in the final policy that solves the assignment problem as an \gls{ip}.

\FloatBarrier
\section{Problem Formulation}\label{sec:formulation}

We consider the problem of controlling a pooling-enabled \gls{rhs} as described in \Cref{sec:problem_statement}. We model vehicle movements on a weighted directed graph $G=(\V,A)$ where $\V$ is a set of nodes and $A$ is a set of arcs. Each node $v \in \V$ represents a location that can act as the origin and destination of requests. Each arc $(u,v) \in A$ represents a drivable road segment from location $u$ to $v$. We assume that the graph is connected. We also assume that the time it takes a vehicle to travel the shortest path from node $u \in V$ to node $v \in V$ is deterministic and is given by $\tau(u,v)$.
Extending our model to the case where travel times are random and vehicles can enter or leave service at random (e.g., due to accidents) is straightforward.

Furthermore, we consider a finite time horizon $T$ with a time-triggered decision epoch $\tInSetTimeSteps$. In every decision epoch $t$, we observe the system \textit{state}, apply a \textit{decision} on each vehicle, receive a \textit{reward}, observe new \textit{exogenous information}, and advance to a new state in the next decision epoch $t+1$ via the \textit{transition functions}.

\subsection{States}\label{subsec:states}
The system state
variable
describes the state of every vehicle in the fleet and the attributes of the travel requests.
The state variable at time $t$ is given by the vector $S_t = (R_t, D_t)$ where
\begin{itemize}
    \item $R_t=(R_{ta})_{a \in \A}$ is the number of vehicles (resources) with attribute $a \in \A$, and
    \item $D_t=(D_{tb})_{b \in \B}$ is the number of travel requests (demand) with attribute $b \in \B$,
\end{itemize}
where the attribute vectors $a$ and $b$ are defined as
\begin{align*}
    a = \begin{pmatrix} o_a \\ d_a \\ l_a \\ n_a \\ \taua \end{pmatrix}
    = 
    \begin{pmatrix} \textnormal{current location} \\ \textnormal{destination} \\ \textnormal{driving range} \\ \textnormal{capacity} \\ \textnormal{actionable time}\end{pmatrix}
\quad\text{and}\quad
    b = \begin{pmatrix} o_b \\ d_b \\ n_b \\ \taur_b \\ \taup_b \\ f_b \end{pmatrix}
    = 
    \begin{pmatrix} \textnormal{origin} \\ \textnormal{destination} \\ \textnormal{headcount} \\ \textnormal{latest response time} \\ \textnormal{latest pickup time} \\ \textnormal{fare} \end{pmatrix}.
\end{align*}
The vehicle's current location $o_a\in \V$ and destination $d_a\in \V$ are nodes in the graph $G$. Let $l^{\max}$ be the maximum driving range. We define $n^{\max}$ to be the maximum capacity of a vehicle, such that $n_a=n^{\max}$ when the vehicle is empty and $n_a=0$ when the vehicle is fully occupied. We distinguish between two types of vehicles in the system:
(i) \emph{empty} vehicles $\Aempt = \{ a \in \A \mid n_a=n^{\max} \wedge o_a = d_a\}$, and (ii) \emph{occupied} vehicles $\Aoccu = \{ a \in \A \mid n_a < n^{\max}\wedge o_a \neq d_a\}$. The actionable time $\taua \in \mathbb{R}_{\geq t}$ is the time at which a vehicle $a$ can start performing a new task. When $\taua>t$, the vehicle $a$ is currently performing a task,
and $\taua=t$ means that the vehicle can start performing a new task immediately.

The trip's origin $o_b\in \V$ and destination $d_b\in \V$ are nodes in the graph $G$. We only consider requests where $o_b \neq d_b$. The headcount $n_b \in \mathbb{Z}_{>0}$ is the number of passengers in request $b$, and we only consider requests where $n_b \leq n^{\max}$. The latest response time $\taur_b \in \{1,\dots,T\}$ is the latest time at which the \gls{rhs} can accept request $b$ by assigning it to a vehicle. The latest pickup time $\taup_b \in \{1,\dots,T\}$ is the latest time at which a vehicle dispatched by the \gls{rhs} can reach the origin of request $b$ and pick up its passengers. At time $t$, the set $\B$ consists only of requests $b$ such that $t \leq \taur_b$ and $t \leq \taup_b$. The trip fare $f_b \in \mathbb{R}_{>0}$ is the amount of money the \gls{rhs} receives if it accepts request $b$.

\subsection{Decisions}\label{subsec:decisions3}
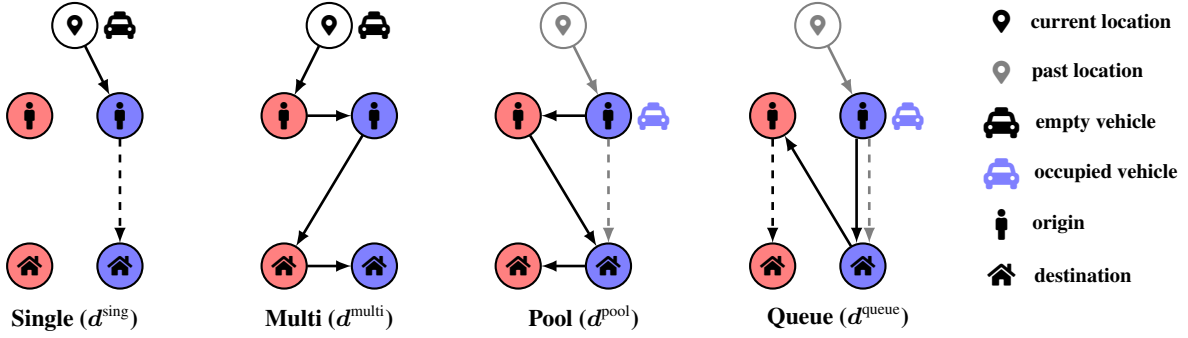
\begin{figure}[htbp]
    \centering
    \tikzset{%
        base/.style args={sized #1}{circle, draw=black, fill=white, thick, inner sep=0pt,minimum size=#1    , align=center},
        base/.default={sized 15mm},
        node/.style={circle, draw=black, fill=white, thick, inner sep=0pt,minimum size=1mm, maximum size=3mm, align=center},
        edge/.style={-latex, line width=1pt, black},
        arc/.style={thick,double arrow=3pt colored by white and black}
    }
    \begin{subfigure}[t]{0.18\textwidth}
        \centering
        \begin{tikzpicture}[scale=0.4, transform shape]
            \node at (5.5,7) {\emptycar};
            \node [base] (a1) at (4,7) {\currentlocation};
            \node [base, fill=blue!50] (o1) at (5.5,+4.0) {\originnodem};
            \node [base, fill=red!50]  (o2) at (2.5,+4.0) {\originnodem};
            \node [base, fill=blue!50] (d1) at (5.5,-1.0) {\destinationnode};
            \node [base, fill=red!50]  (d2) at (2.5,-1.0) {\destinationnode};
            \draw[edge, solid] (a1) -- (o1);
            \draw[edge, dashed] (o1) -- (d1);
        \end{tikzpicture}
        \caption*{\textbf{Single ($\bm{d^{\textnormal{sing}}}$)} 
        }
    \end{subfigure}
    \hfill
    \begin{subfigure}[t]{0.18\textwidth}
        \centering
        \begin{tikzpicture}[scale=0.4, transform shape]
            \node at (5.5,7) {\emptycar};
            \node [base] (a1) at (4,7) {\currentlocation};
            \node [base, fill=blue!50] (o1) at (5.5,+4.0) {\originnodem};
            \node [base, fill=red!50]  (o2) at (2.5,+4.0) {\originnodem};
            \node [base, fill=blue!50] (d1) at (5.5,-1.0) {\destinationnode};
            \node [base, fill=red!50]  (d2) at (2.5,-1.0) {\destinationnode};
            \draw[edge, solid] (a1) -- (o2);
            \draw[edge, solid] (o2) -- (o1);
            \draw[edge, solid] (o1) -- (d2);
            \draw[edge, solid] (d2) -- (d1);
        \end{tikzpicture}
        \caption*{\textbf{Multi ($\bm{d^{\textnormal{multi}}}$)} 
        }
    \end{subfigure}
    \hfill
    \begin{subfigure}[t]{0.18\textwidth}
        \centering
        \begin{tikzpicture}[scale=0.4, transform shape]
            \node at (7,4) {\occupiedcar};
            \node [base, draw=gray] (a1) at (4,7) {\pastlocation};
            \node [base, fill=blue!50] (o1) at (5.5,+4.0) {\originnodem};
            \node [base, fill=red!50]  (o2) at (2.5,+4.0) {\originnodem};
            \node [base, fill=blue!50] (d1) at (5.5,-1.0) {\destinationnode};
            \node [base, fill=red!50]  (d2) at (2.5,-1.0) {\destinationnode};
            \draw[edge, solid, gray] (a1) -- (o1);
            \draw[edge, dashed, gray] (o1) -- (d1);
            \draw[edge, solid] (o1) -- (o2);
            \draw[edge, solid] (o2) -- (d1);
            \draw[edge, solid] (d1) -- (d2);
        \end{tikzpicture}
        \caption*{\textbf{Pool ($\bm{d^{\textnormal{pool}}}$)}}
    \end{subfigure}
    \hfill
    \begin{subfigure}[t]{0.18\textwidth}
        \centering
        \begin{tikzpicture}[scale=0.4, transform shape]
            \node at (7,4) {\occupiedcar};
            \node [base, draw=gray] (a1) at (4,7) {\pastlocation};
            \node [base, fill=blue!50] (o1) at (5.5,+4.0) {\originnodem};
            \node [base, fill=red!50]  (o2) at (2.5,+4.0) {\originnodem};
            \node [base, fill=blue!50] (d1) at (5.5,-1.0) {\destinationnode};
            \node [base, fill=red!50]  (d2) at (2.5,-1.0) {\destinationnode};
            
            \draw[edge, solid, gray] (a1) -- (o1);
            \draw[edge, dashed, gray] ([xshift=6pt]o1.south) -- ([xshift=6pt]d1.north);
            \draw[edge, solid] ([xshift=-6pt]o1.south) -- ([xshift=-6pt]d1.north);

            \draw[edge, solid] (d1) -- (o2);
            \draw[edge, solid, dashed] (o2) -- (d2);
        \end{tikzpicture}
        \caption*{\textbf{Queue ($\bm{d^{\textnormal{queue}}}$)}}
    \end{subfigure}
    \begin{subfigure}[t]{0.2\textwidth}
        \centering
        \begin{tikzpicture}[scale=0.45, transform shape, yshift=-5cm]
            \node at (0,7.5) {\currentlocation};
            \node at (0,6) {\pastlocation};
            \node at (0,4.5) {\emptycar};
            \node at (0,3) {\occupiedcar};
            \node at (0,1.5) {\originnodem};
            \node at (0,0) {\destinationnode};
            
            \node[align=left] at (2.7,7.5) {{\LARGE \textbf{ \; current location}}};
            \node[align=left] at (2.3,6) {{\LARGE \textbf{ \; past location}}};
            \node[align=left] at (2.55,4.5) {{\LARGE \textbf{ \; empty vehicle}}};
            \node[align=left] at (2.85,3) {{\LARGE \textbf{ \; occupied vehicle}}};
            \node[align=left] at (1.45,1.5) {{\LARGE \textbf{ \; origin}}};
            \node[align=left] at (1.85,0) {{\LARGE \textbf{ \; \; \; destination}}};
        \end{tikzpicture}
    \end{subfigure}
    
    \caption{Four types of decisions that assign a vehicle to serve one or more requests. 
        Solid arcs must be traversed and dashed arcs are tentative routes that may be modified later.
        Black arcs are part of the post-decision state.
        Gray dashed arcs are part of the pre-decision state, and gray solid arcs are part of the path history.
    }
    \label{fig:decisions}
\end{figure}

A decision $d \in \D$ is a task that a vehicle can be assigned to perform.
For all vehicles with attribute $a\in\A$ and requests with attribute $b\in\B$, let $\D(a)$ be the subset of decisions that can be assigned to a vehicle $a \in \A$. 
Let $\D(a,b) \subseteq \D(a)$ be the set of decisions that assign vehicle $a$ to serve request $b$, see Figure \ref{fig:decisions} for an illustration.
The set of all decisions is given by $\D = \bigcup_{a \in \A} \D(a)$, where
\begin{align*}
    \D(a) = &\Decs_a \cup \Decm_a \cup \Decp_a \cup \Decq_a \cup \Decrel_a \cup \Deccont_a \cup \Decidle_a \cup \Decrec_a.
\end{align*}
To define the individual subsets, we introduce the following notation. Let $P$ be the set of paths in $G$ and let $\sigma$ be a function that maps each path $p\in P$ to the driving range that a vehicle requires to traverse $p$. A path $p\in P$ is $a$-\emph{feasible} for vehicle $a\in \mathcal{A}$ if it starts with $o_a$ and $l_a \geq \sigma(p)$. Let $P_a$ be the set of $a$-feasible paths, for each $a\in \mathcal{A}$, and let  $B\subseteq\mathcal{B}$ be a nonempty set of requests. A path $p\in P_a$ that $a\in\A$ traverses to satisfy the requests in $b \in B$ is $a$-$B$-\emph{feasible} if it satisfies the following conditions: (1) for every request $b \in B$, $p$ visits the origin $o_b$ before the destination $d_b$, and $p$ reaches the origin $o_b$ before the latest pickup time $\taup_b$; and (2) vehicle $a$ has enough seats to pick up and drop off the requests $B$ in the order of $p$.
For each $a\in\mathcal{A}$ and $B\subseteq\mathcal{B}$, let $P_{aB}$ the set of $a$-$B$-feasible paths and let $P_{aB}(u,v) \in P_{aB}$ be the set of $a$-$B$-feasible paths that visit location $u\in\V$, then location $v\in\V$, and may or may not visit other locations before $u$, after $v$, and in between.

\paragraph{Single trip:} We use $\decs_b$ to denote the decision to assign an empty vehicle to serve a single travel request $b \in \B$.
We define
\begin{equation}\label{eq:Decs}
    \Decs_a = \left\{\decs_b \:\Big\vert\: b \in \B \wedge P_{a\{b\}}\neq\emptyset\right\}
\end{equation}
if $a \in \Aempt$. Otherwise, $\Decs_a = \emptyset$.

\paragraph{Multi-trip:} We use $\decm_{Bp}$ to denote the decision to assign an empty vehicle to serve a set of requests $B \subseteq \B$ on the path $p\in P_{aB}$. We define
\begin{equation}\label{eq:Decm}
    \Decm_a = \left\{\decm_{Bp} \:\Big\vert\: B \subseteq \B  \wedge |B| > 1 \wedge p\in P_{aB}\right\}
\end{equation}
if $a \in \Aempt$. Otherwise, $\Decm_a = \emptyset$.

\paragraph{Pool:} We use $\decp_{bp}$ to denote the decision to assign an occupied vehicle $a\in\mathcal{A}$ to serve a new travel request $b\in \B$, together with the passengers that are already in the vehicle, on the path $p\in P_{a\{b\}}$.
We define
\begin{equation}\label{eq:Decp}
    \Decp_a = \left\{\decp_{bp} \:\Big\vert\: b \in \B \wedge n_a+n_b \leq n^\text{max} \wedge p \in P_{a\{b\}}(o_b, d_a) \right\}
\end{equation}
if $a \in \Aoccu$. Otherwise, $\Decp_a = \emptyset$.

\paragraph{Queue:} We use $\decq_b$ to denote the decision to assign an occupied vehicle to serve a new travel request $b \in \B$ after it drops off its current passengers. We define
\begin{equation}\label{eq:Decq}
    \Decq_a = \left\{\decq_b \mid b \in \B \wedge P_{a\{b\}}(d_a, o_b)\neq\emptyset \right\}
\end{equation}
if $a \in \Aoccu$. Otherwise, $\Decq_a = \emptyset$.

\paragraph{Relocate:} We use $\decrel_v$ to denote the decision to relocate a vehicle at location $v \in \V$. Vehicles can either relocate to adjacent locations or locations that are reachable before the next decision epoch. Vehicles $a\in\A$ whose actionable time $t_a$ is after the next decision epoch $t$ cannot relocate because the decision would be premature. We define
\begin{equation}\label{eq:Decrel}
    \Decrel_a = \left\{\decrel_v \:\Big\vert\: v \in \V  \wedge \tau(o_a,v) \leq l_a \wedge t_a < t+1 \\ \wedge \big((o_a, v) \in A \vee \tau(o_a,v) \leq t+1\big) \right\}
\end{equation}
if $a \in \Aempt$. Otherwise, $\Decrel_a = \emptyset$.

\paragraph{Continue:} We use $\deccont$ to denote the decision to let an already occupied vehicle continue its task with no new instructions. We define $\Deccont_a =\{\deccont\}$ if $a\in\Aoccu$, and $\Deccont_a =\emptyset$ otherwise.

\paragraph{Idle:} We use $\decidle$ to denote the decision to park an empty vehicle. We define $\Decidle_a = \{\decidle\}$ if $a \in \Aempt$ at time $\taua$. Otherwise, $\Decidle_a = \emptyset$.

\paragraph{Recharge:} We use $\decrec$ to denote the decision to refuel or recharge a vehicle. Vehicles $a\in\A$ whose actionable time $t_a$ is after the next decision epoch $t$ cannot recharge because the decision would be premature. We define 
\begin{equation}\label{eq:Decrec}
    \Decrec_a = \{\decrec \mid l_a < l^{\max} \wedge t_a < t+1\} 
\end{equation}
if $a \in \Aempt$. Otherwise, $\Decrec_a = \emptyset$.


We define a decision variable $x_{tad}$ as the number of times we apply a decision $d \in \D(a)$ to a vehicle $a \in \A$ at time $t$. The decision vector $x_t = (x_{tad})_{a \in \A, d \in \D(a)}$ must respect the constraints
\begin{subequations}
\begin{align}
    \sum_{d \in \D(a)} x_{tad} &= R_{ta} && \forall a \in \A, \label{eq:flow_balance} \\
    \sum_{a\in\A} \sum_{d \in \D(a,b)} x_{tad} &\leq D_{tb} && \forall b \in \B, \label{eq:feasible_assignment} \\
    x_{tad} &\in \mathbb{Z}_{\geq 0} && \forall a \in \A,\ d \in \D(a). \label{eq:integrality}
\end{align}
\end{subequations}
Constraint \eqref{eq:flow_balance} is a flow-balance constraint. Constraint \eqref{eq:feasible_assignment} ensures that the number of vehicles assigned to a travel request with attribute $b \in \B$ does not exceed the number of available requests with such attribute. Constraint \eqref{eq:integrality} ensures that the number of vehicles $a \in \A$ assigned to a decision $d \in \D(a)$ is non-negative and integer. Constraints \eqref{eq:flow_balance}, \eqref{eq:feasible_assignment}, and \eqref{eq:integrality} define the feasible set $\X(S_t)$.

\subsection{Rewards}\label{subsec:rewards}

The reward is given by the function
\begin{equation}\label{eq:reward_joint}
    C(x_t) = \sum_{a \in \A}\sum_{d \in \D(a)} c_{tad} x_{tad},
\end{equation}
where the contribution $c_{tad}$ is defined as follows.
Assigning a vehicle $a\in\A$ to serve a request $b\in \B$ via a decision $d \in \{\decs_b, \decq_b\}$ at time $t$ produces a contribution that is equal to the trip fare, i.e., $c_{tad} = f_b$.
Let $\detourPenalty(a,d)\in\mathbb{R}_{\leq 0}$ be the detour penalty incurred when a vehicle $a$ is assigned to decision $d \in \D(a)$. The penalty is equal to zero if $d$ does not involve any detour.
Assigning a vehicle $a\in\A$ to serve a request $b\in \B$ on path $p\in P_{a\{b\}}$ via decision $d = \decp_{bp}$ produces the contribution $c_{tad} = f_b + \detourPenalty(a,d)$.
Assigning a vehicle $a\in\A$ to serve a set of requests $B\in \B$ on a path $ p\in P_{aB}$ via decision $d = \decm_{Bp}$ produces the contribution 
\begin{equation}\label{eq:multi_trip_contribution}
    c_{tad} = \sum_{b\in B} f_b + \detourPenalty(a,d) 
\end{equation}
Let $\rechargeCost(a)\in\mathbb{R}_{\leq 0}$ be the cost of fully recharging a vehicle $a \in \A$, which depends on the vehicle's current range $l_a$. Assigning a vehicle $a\in\A$ to a recharge decision $d \in \Decrec_a$ at time $t$ produces the contribution $c_{tad} = \rechargeCost(a)$. 
Assigning a vehicle $a\in\A$ to an idle, relocate, or continue decision $d \in \Decidle_a \cup \Decrel_a \cup \Deccont_a$ at time $t$ produces no contribution or penalty, i.e., $c_{tad} = 0$.

\subsection{Exogenous information}\label{subsec:exagogenous_information}
The exogenous information is the new information that the \gls{rhs} receives between decision epochs $t-1$ and $t$. The exogenous information is given by the demand vector $\hat D_t=(\hat D_{tb})_{b \in \B}$ where $\hat D_{tb}$ is the number of travel requests with attribute $b \in \B$ that first became known between time $t-1$ and $t$.

\subsection{Transition functions}\label{subsec:transition_functions}
The transition functions specify how the system evolves from one state to the next.
Let $a' = a^{M}(a,d)$ be the future attribute of a vehicle $a$ to which a decision $d\in\D(a)$ is applied. The function $a^{M}$ defines how the attribute of vehicle $a$ changes after it is assigned to decision $d$.
We define a \emph{post-decision state} variable $S_t^x$. The post-decision state $S_t^x$ is the state of the system after we make decision $x_t$ but before any new information (i.e., exogenous demand $\hat D_{t+1}$) has arrived. 
The post-decision state is defined as $S_t^x = (R_t^x, D_t^x)$ where
\begin{itemize}
    \item $R^x_t=(R^x_{ta'})_{a' \in \A}$ and $R^x_{ta'} = \sum_{a \in \A}\sum_{d \in \D(a)} x_{tad} \mathbf{1}_{\{a^{M}(a,d)=a'\}}$.
    \item $D^x_t=(D^x_{tb})_{b \in \B}$ and $D^x_{tb} = \left(D_{tb} - \sum_{a \in \A} \sum_{d_b \in \D(a,b)} x_{tad_b}\right) \mathbf{1}_{\{\taur_b \geq t+1\}}$.
\end{itemize}
We can obtain $D^{x}_{tb}$ from the vector of slack variables of the demand constraints \eqref{eq:feasible_assignment} by setting every element that corresponds to a request $b\in\B$ with $\taur_b \leq t+1$ to 0.
We now complete the definition of the transition function as
\begin{align*}
    S_{t+1} &= S^{M}(S_t,x_t,\hat D_{t+1}) = (R^x_{t}, D^x_{t}+\hat D_{t+1}) = (R_{t+1}, D_{t+1})
\end{align*}
where $R_{t+1}$ and $D_{t+1}$ are the \textit{resource} and \textit{demand} vectors at time $t+1$, respectively. We refer the reader to \Cref{sec:transition_function_details} for a detailed definition of the function $a^{M}(a,d)$.



\section{Policies}\label{sec:policies}
A policy is a function $\pi: S_t \to \mathcal{X}(S_t)$ that maps the state $S_t$ to a decision vector $x_t \in \mathcal{X}(S_t)$. Let $\Pi$ be the set of policies. The optimal policy $\pi^*$ maximizes the sum of expected rewards over the entire planning horizon, i.e.,
\begin{equation}\label{eq:optimal_policy_definition}
    \pi^* = \argmax_{\pi \in \Pi} \sum_{t=1}^{T} \mathbb{E}\left[\left\{ C(\pi(S_t)) \; | \; S_0 \right\}\right].
\end{equation}

\subsection{Myopic policy}\label{subsec:myopic_policies}
A myopic policy makes decisions by maximizing the reward function~\eqref{eq:reward_joint} subject to the constraints~\eqref{eq:flow_balance}--\eqref{eq:integrality}, i.e.,
\begin{equation}\label{eq:myopic_policy}
    \myopicpolicy{}(S_t) = \argmax_{x_t \in \X(S_t)}\ \sum_{a \in \A}\sum_{d \in \D(a)} c_{tad} x_{tad}.
\end{equation}
The myopic policy in~\eqref{eq:myopic_policy} is suboptimal because it does not consider how a decision $x_t$ affects the next state $S_{t+1}$ the system will evolve into and, through that, all subsequent rewards. If the fleet of vehicles is electric, a simple way to improve $\myopicpolicy{}$ is to introduce a threshold (parameter) $\theta \in (0,1]$ that forces the vehicle to charge only when its battery level is below the threshold. This can be done by modifying Problem~\eqref{eq:myopic_policy} to become
\begin{equation}\label{eq:param_myopic_policy}
\begin{aligned}
\pi^\text{PM}(S_t) = \argmax_{x_t \in \X(S_t)} \quad 
& \sum_{a \in \A} \sum_{d \in \D(a)} c_{tad} x_{tad} \\
\text{s.t.} \quad 
& x_{tad} = R_{ta} \cdot \mathbf{1}_{\{l_a < \theta l^{\max}\}}, \quad \forall\, a \in \A,\; d \in \Decrec(a).
\end{aligned}
\end{equation}
We refer to the policy in~\eqref{eq:param_myopic_policy} as the \gls{pm} policy.

\subsection{Value function approximation policy}
\label{subsec:vfa_policies}
To capture the downstream impact of current decisions on future rewards, we use the well-known Bellman equation. We define the value function $V(S_t)$ as the expected total reward from time $t$ to $T$ given that the system is in state $S_t$ at time $t$, i.e.,
\begin{equation}\label{eq:bellman_equation}
    V(S_t) = \max_{x_t \in \X(S_t)} \big\{ C(x_t) + \mathbb{E}[V(S_{t+1}) \mid S_t, x_t] \big\},
\end{equation}
for $\tInSetTimeSteps$ and $V(S_{T+1})=0$.
A solution to the Bellman equation \eqref{eq:bellman_equation} gives the optimal policy in \eqref{eq:optimal_policy_definition}.
Solving the Bellman equation~\eqref{eq:bellman_equation} exactly is intractable. A common approach to address this challenge is to approximate the value function $V(S_t)$ by a function 
that is easier to compute.
We use a linear approximation \citep{simao:2009} given by
\begin{equation}\label{eq:linear_approximation}
    \bar{V}^x(S^{x}_t) = \sum_{a \in \A} \bar{v}^{R}_{a}R^{x}_{ta} + \sum_{b\in\B} \bar{v}^{D}_{tb}D^{x}_{tb}
\end{equation}
where 
\begin{itemize}
    \item $\bar{v}^{R}_{a}= \frac{\partial V(S_t)}{\partial R_{ta}}$ is the marginal value of a resource with attribute $a\in \A$ at any time $\tInSetTimeSteps$, and 
    \item $\bar{v}^{D}_{tb}= \frac{\partial V(S_t)}{\partial D_{tb}}$ the marginal value of a demand with attribute $b\in \B$ at time $t$.
\end{itemize}
Note that we omit the time index $t$ in $\bar{v}^{R}_{a}$ since the vehicle's actionable time is already encoded in the attribute $a$.

Substituting $\bar{V}^x(S^{x}_t)$ into \Cref{eq:bellman_equation}, we obtain the \gls{vfa} policy 
\begin{align}\label{eq:vfa_policy}
        \pi^\mathrm{VFA}(S_t) = \argmax_{x_t \in \X(S_t)}
        \left\{ 
             \sum_{a \in \A}\sum_{d \in \D(a)}
             \left(c_{tad} +\bar{v}^{R}_{a^{M}(a,d)} \right)x_{tad}
            +\sum_{b\in\B} \bar{v}^{D}_{tb}D^{x}_{tb}
        \right\}.
\end{align}
We refer the reader to \Cref{sec:policies_details} for a detailed derivation of \Cref{eq:vfa_policy}.

\subsection{Policies design and implementation}
\label{subsec:policies_design_implementation}
Solving the \gls{pm} problem in \eqref{eq:myopic_policy} or the \gls{vfa} problem in \eqref{eq:vfa_policy} is associated with several challenges that require careful design decisions. We discuss and address these challenges in the remainder of this section.
\subsubsection{Enumerating the multi-trip decisions} 
The presence of the multi-trip decision $\decm$ leads to an exponential increase in the complexity of enumerating the elements of the decision set $\D(a)$ for each $a \in \A$. Not only do we need to enumerate all possible subsets of requests $B \subseteq \B$ that can be assigned to vehicle $a$, but we also need to consider the different feasible paths that a vehicle with attribute $a$ can take to serve the requests in $B$.

To enumerate the subsets of requests that a vehicle $a\in \A$ can serve, we use a modified version of the efficient \emph{RTV graph algorithm} of \cite{alonso:2017}. The original RTV graph algorithm creates subsets of requests $B \subseteq \B$ that can be served by a virtual empty vehicle located at the same location as the origin of one of the requests in $B$, while satisfying all operational constraints (e.g., time windows, vehicle capacity, maximum waiting time). To limit computation, trips are constructed incrementally, starting with single-request trips and adding one request at a time until no more requests can be added without violating any operational constraint or exceeding a pre-specified maximum number of requests per trip. Then, the algorithm discards each subset of requests $B$, for which there exists no vehicle that can reach the origin of at least one request in $B$ before its latest pickup time. 

We modify the RTV graph algorithm by adding all constraints of our model.
Our modified algorithm enumerates all feasible trips $B \in \B$ separately for each vehicle attribute $a \in \A$, again proceeding incrementally. Because of the range constraint, checking whether vehicle $a$ can serve trip $B$ requires solving a \gls{spprc} \citep{Irnich2005}. Solving this \gls{spprc} yields a shortest path $p^* \in P_{aB}$ that serves all requests in $B$ while satisfying all operational constraints. 
This problem can be solved very efficiently using a labelling algorithm since we have no negative-weight arcs. 
If we assume that the detour penalty function $\detourPenalty{}$ is monotonically increasing in the path length, we can limit the number of paths $|P_{aB}|$ that we need to consider for each subset of requests $B$. We formalize this assumption in the following definition.
\detourPenaltyMonoDefinition{}
Since the \gls{pm} policy does not consider future states, we only need to include one multi-trip decision $\decm_{Bp^*}$ per feasible subset of requests $B$ in the decision set $\D(a)$ for each vehicle $a \in \A$. In this case, $p^* \in P_{aB}$ is one of the shortest paths in $P_{aB}$.
That is, since the \gls{pm} policy is indifferent to the location, battery level, and time of arrival of vehicle $a$ after serving the requests in $B$, we only need to include the decision that maximizes the immediate contribution $c_{tad}$ given by \Cref{eq:multi_trip_contribution}. 
Assuming that the detour penalty function $\detourPenalty{}$ is monotonically increasing in the path length, the shortest path $p^* \in P_{aB}$ will yield the smallest detour penalty $\detourPenalty(a,\decm_{Bp^*})$ among all paths in $P_{aB}$, which in turn maximizes the contribution $c_{tad}$. This shortest path $p^* \in P_{aB}$ is the solution of the \gls{spprc} problem that we solve in our modified RTV graph algorithm.

The situation for the \gls{vfa} policy is more involved because \eqref{eq:vfa_policy} maximizes not only the immediate contribution $c_{tad}$ but also the downstream value captured by the term $\bar{v}^{R}_{a^{M}(a,d)}$.
For example, a path $p'$ that ends at an airport with high demand may be preferred over a shorter $p^*$ that ends in a suburban neighborhood with little demand, even if $p'$ takes longer to complete and yields a lower immediate contribution $c_{tad}$ than $p^*$.
We show that enumerating all possible paths $p \in P_{aB}$ for each subset of requests $B \subseteq \B$ in the decision set $\D(a)$ is not necessary to obtain an optimal solution to \eqref{eq:vfa_policy}. Instead, for each subset of requests $B \subseteq \B$, it suffices to consider at most $|B|$ paths from $P_{aB}$. Specifically, for each subset of requests $B$, we need to consider at most one path $p^*_b \in P_{aB}$ per each request $b \in B$, which is one of the shortest paths in $P_{aB}$ that end at the request destination $d_b$. Fortunately, without requiring any additional enumeration, the labeling algorithm we use to solve the \gls{spprc} yields a set $\{p^*_b \in P_{aB} \, \mid \, b \in B\}$, which contains all paths that we need to consider. We formalize this result in the following definitions and propositions.
\enumerationDefinitions{}
\partialOrderDefinition{}
\Cref{def:partial_order} allows us to adopt and extend the monotonicity results of \cite{alkanj:2020} to our setting in the following proposition. 
\monotonicityProposition{}
\begin{proof}
    See \Cref{sec:policies_details}.
\end{proof}
\enumerationProposition{}
\begin{proof}
    See \Cref{sec:proof}.
\end{proof}

\subsubsection{Obtaining the marginal attribute values}    
\label{sec:adp}
To implement the \gls{vfa} policy in \eqref{eq:vfa_policy}, we need to know the marginal attribute values $\bar{v}^{R}_{a}$ for every $a \in \A$ and $\bar{v}^{D}_{tb}$ for every $b \in \B$ and $\tInSetTimeSteps$. In practice, these values are unknown and need to be estimated.
We estimate them using a Forward \gls{adp} algorithm that starts with initial values $\bar{v}^{R,0}_{a}$ and $\bar{v}^{D,0}_{tb}$ and iteratively improves these estimates by simulating the \gls{rhs} for $N$ iterations. 
In the $n$-th iteration, we obtain the decision $x^n_t$ and the estimates $\hat{v}^{R,n}_{ta}$ and $\hat{v}^{D,n}_{tb}$ of the slope of the objective function with respect to $R_{ta}$ and $D_{tb}$, respectively, at the points $R^x_{ta}$ and $D^x_{tb}$, respectively. We then update the current marginal attribute value estimates $\bar{v}^{R,n-1}_a$ and $\bar{v}^{D,n-1}_{tb}$ to obtain new estimates $\bar{v}^{R,n}_a$ and $\bar{v}^{D,n}_{tb}$. The last step of the iteration is to transition to the next state using $x^n_t$ and a random demand realization. As the final result of the algorithm, we obtain the values $\bar{v}^{R,N}_{a}$ for every $a\in\A$ and $\bar{v}^{D,N}_{tb}$ for every $b\in\B$ and $\tInSetTimeSteps$. We then use these values in \eqref{eq:vfa_policy} to obtain our \gls{vfa} policy. 
Note that we retain the time index $t$ in $\hat{v}^{R,n}_{ta}$ to emphasize that it is an estimate of the objective function's slope at the point $R^x_{ta}$ at time $t$ in iteration $n$.
\Cref{sec:policies_details} gives a detailed description of the algorithm.

To use the algorithm, we need a surrogate policy that allows us to obtain the estimates $\hat{v}^{R,n}_{ta}$ and $\hat{v}^{D,n}_{tb}$ in the $n$-th iteration. 
If the feasible region $\X$ of \eqref{eq:vfa_policy} has a totally unimodular coefficient matrix and an integer right-hand side, we can relax the integrality constraints \eqref{eq:integrality} and solve \eqref{eq:vfa_policy} as an \gls{lp}. Doing so allows us to obtain the marginal values $\hat{v}^{R,n}_{ta}$ and $\hat{v}^{D,n}_{tb}$ directly from the dual variables associated with the resource constraints \eqref{eq:flow_balance} and demand constraints \eqref{eq:feasible_assignment}, respectively.
To obtain a feasible region $\X^\mathrm{LP}$ that has a totally unimodular coefficient matrix, we remove the multi-trip decision $\decm$ from $\X$. We only use $\X^\mathrm{LP}$ in the surrogate policy, and we only use the surrogate policy to obtain $\bar{v}^{R,N}_{a}$ and $\bar{v}^{D,N}_{tb}$. Once we have these values, we can use the final \gls{vfa} policy that solves \eqref{eq:vfa_policy} with the multi-trip decisions. That is, we include the multi-trip decisions in the final \gls{vfa} policy but not in the surrogate policy used in the \gls{adp} procedure.

\subsubsection{Exploitation versus exploration trade-off}
\label{sec:exploration_exploitation}
\gls{adp} faces the classic \textit{exploration vs.\ exploitation trade-off}.
On the one hand, we want to \textit{explore} different decisions to visit a diverse set of states and learn their values. On the other hand, we want to \textit{exploit} the current estimates of the value function to make meaningful progress in the learning process.
For instance, without exploration, the policy may get stuck in a local optimum that allocates all vehicles to the same location repeatedly because it only knows the values of that location and never visits other locations.
Conversely, without exploitation, the policy may wander randomly across the state space without making meaningful progress in learning the values of the states.

A standard approach to address this trade-off is to introduce some randomness in the decision-making process to encourage exploration.
We follow this approach by adding some randomness to the relocate decisions made by the \gls{vfa} policy in \eqref{eq:vfa_policy}. Specifically, after solving Problem~\eqref{eq:vfa_policy} in each decision epoch $t$, we post-process the chosen relocate destinations, as follows. We consider the subset of attributes for which at least one relocate decision variable is non-zero, i.e., $\A^\mathrm{relc} = \{a\in\A\ \mid r_a \geq 0 \},$ where $r_a = \sum_{d\in\Decrel_a} x_{tad}$.
For each attribute $a\in\A^\mathrm{relc}$, we draw, with replacement, $r_a$ random samples from the set of feasible relocate destinations of $a$. In the $n$-th iteration, the probability of drawing a destination $v\in\V$ is proportional to the value $\bar{v}_{a'}^{R,n}$ of the post-decision state $a' = a^M(a,\decrel_v)$.





Another approach to address the exploration vs.\ exploitation trade-off is to leverage the monotonicity properties introduced in \Cref{prop:monotonicity} to update some of the $\bar{v}^{R,n}$ values during the \gls{adp} procedure. 
For instance, suppose we find that our current estimate of the marginal attribute value $\bar{v}^{R,n}_{a} > \bar{v}^{R}_{a'}$ for two attributes $a, a' \in \A$, even though $a \preccurlyeq a'$. If we leave theses estimates as-is, the policy may wrongly exploit the supposedly higher value of attribute $a$ by more frequently choosing decisions whose post-decision state is $a$. We can correct this inconsistency by setting $\bar{v}^{R,n}_{ta'} = \bar{v}^{R,n}_{ta}$. This correction encourages exploration by increasing the value of attribute $a'$, which may lead the policy to more frequently choose decisions whose post-decision state is $a'$, thereby exploring states associated with $a'$.

\FloatBarrier
\section{Numerical Experiments}\label{sec:experiments}
In this section, we provide numerical results and analyze the impact of different features, decisions, and constraints on the performance of the \gls{rhs}. We start by describing the instances we use (\Cref{sec:instances}) and the problem settings (\Cref{sec:settings}). We then present the results of our numerical experiments (\Cref{sec:results}). We refer the reader to \Cref{sec:experiments_details} for implementation details.

\subsection{Benchmark instances}
\label{sec:instances}

Benchmark instances from the literature are either unpublished or do not cover the full range of problem settings that we consider. Furthermore, most existing benchmarks consider only one area and demand distribution. \citet{ulmer:2019} and \citet{yu:2020} show that spatial differences in the demand distribution induce differences in the relative performance of policies. A comprehensive numerical study should, therefore, include multiple areas and demand distributions. We create a set of instances. Upon publication of this paper, we will publish the instances, along with the code we used to create them, as readily usable files on Zenodo.

A popular choice for the street network used in the literature is \gls{NYC}. We create four street networks using data published by the \gls{TLC}. The data comprises real taxi trips from the years 2009 and 2018. Two of our street networks represent the Manhattan borough alone, while the other two represent the four eastern boroughs of \gls{NYC}, which are Manhattan, The Bronx, Brooklyn, and Queens. We use a fine-grained representation of the street networks with the 2009 \gls{TLC} data, and a coarse-grained representation with the 2018 \gls{TLC} data. We associate each street network with a demand distribution. Each demand distribution is represented by a set of sample paths, where each sample path consists of a sequence of travel requests that occur between time $t=0$ and $t=T$. We refer to every combination of one of the four street networks and its corresponding demand distribution as an instance. \Cref{tab:instances} summarizes the characteristics of each instance. We refer the reader to \Cref{sec:street_networks} for a detailed description of how we create the instances.

\begin{table}[htbp]
    \centering
    \resizebox{\textwidth}{!}{%
    \ra{1}
    \begin{tabular}{@{}cclrrrr@{}}
    \toprule
    Instance &  \gls{TLC} data year & Boroughs & $|\V|$ & $|A|$ & \# of Vehicles & Avg.~\# of Requests (per sample path)\\
    \midrule
    1 & 2009 & Manhattan        &  923 & 3581 & 6400 & 395605 \\
    2 &      & four eastern &  9017 & 37511 & 8000 & 424636 \\
    3 & 2018 & Manhattan        &  64 & 319 & 4400 & 240536 \\
    4 &      & four eastern &  231 & 1280 & 6250 & 282739 \\ \bottomrule
    \end{tabular}
    }
    \caption{Characteristics of the benchmark instances.}\label{tab:instances}
\end{table}

\subsection{Problem settings and parameters}
\label{sec:settings}

We consider six different problem settings per instance in a full factorial design. Factors are whether the system is pooling-enabled (i.e., whether $\decp$ and $\decm$ are allowed), and whether the fleet consists of \gls{ICE} vehicles, \glspl{EV} that use \gls{DCFC}, or \glspl{EV} that use slower \gls{l2c}.

We use the following vehicle parameters.
All types of vehicles can carry up to four passengers at a time.
The time it takes to fill the tank or charge the battery (i.e., execute $\decrec$) is
   $\ell(l_a) = (l^{\max} - l_a)\symbolRechargingRate + \qty{15}{min}.$
The slope $\symbolRechargingRate$ depends on the type of vehicle, whereas the intercept is independent of the type of vehicle and accounts for driving to the nearest gas station or charger, queuing, etc.
Fully-fueled \gls{ICE} vehicles can drive for $l^{\max} = \rangeICE$ and completely filling their tank requires \qty{1}{\minute}, i.e., $\symbolRechargingRate = \rechargingRateICE{}$. For \glspl{EV}, we use the \gls{EV} parameters from \citet{alkanj:2020} and convert them to durations, assuming that vehicles either stay at rest or move at the average \gls{NYC} taxi speed, \qty[per-mode=symbol]{18.2}{\kilo\meter\per\hour} according to \gls{TLC} data. We use this constant speed only for determining $\symbolRechargingRate$, not for simulating driving durations (cf. \Cref{sec:experiments_details}). This results in a maximum range of $l^{\max} = \rangeEV$ and a recharging rate of $\symbolRechargingRate{} = \rechargingRateDCFC{}$ for \gls{DCFC} vehicles. We assume that \gls{l2c} vehicles have the same range but recharge 20 times slower, i.e., $\symbolRechargingRate{} = \rechargingRateLTwoC{}$.

We implement one \gls{pm} policy and one \gls{vfa} policy. For each instance, we create 3000 training and 30 test sample paths, along with an initial state for each path. To enable a fair comparison, we reuse the same sample paths and initial states across policies and problem settings.

We use a planning horizon of $T = \qty{24}{h}$, where the \gls{rhs} operator makes a decision every two minutes. This means that the planning horizon consists of 720 decision epochs. We assume energy prices and detour penalties are constant and negligible, i.e., $\rechargeCost = \detourPenalty{} = \text{\$}0$. We consider multi-trip decisions up to $|B| \leq 2$. This is inspired by \citet{alonso:2017}, who report that the average number of passengers per vehicle rarely exceeds two and increased seat capacities yield diminishing returns.

\subsection{Results}\label{sec:results}
\begin{figure}[htbp]  
    \includegraphics[width=0.95\textwidth]{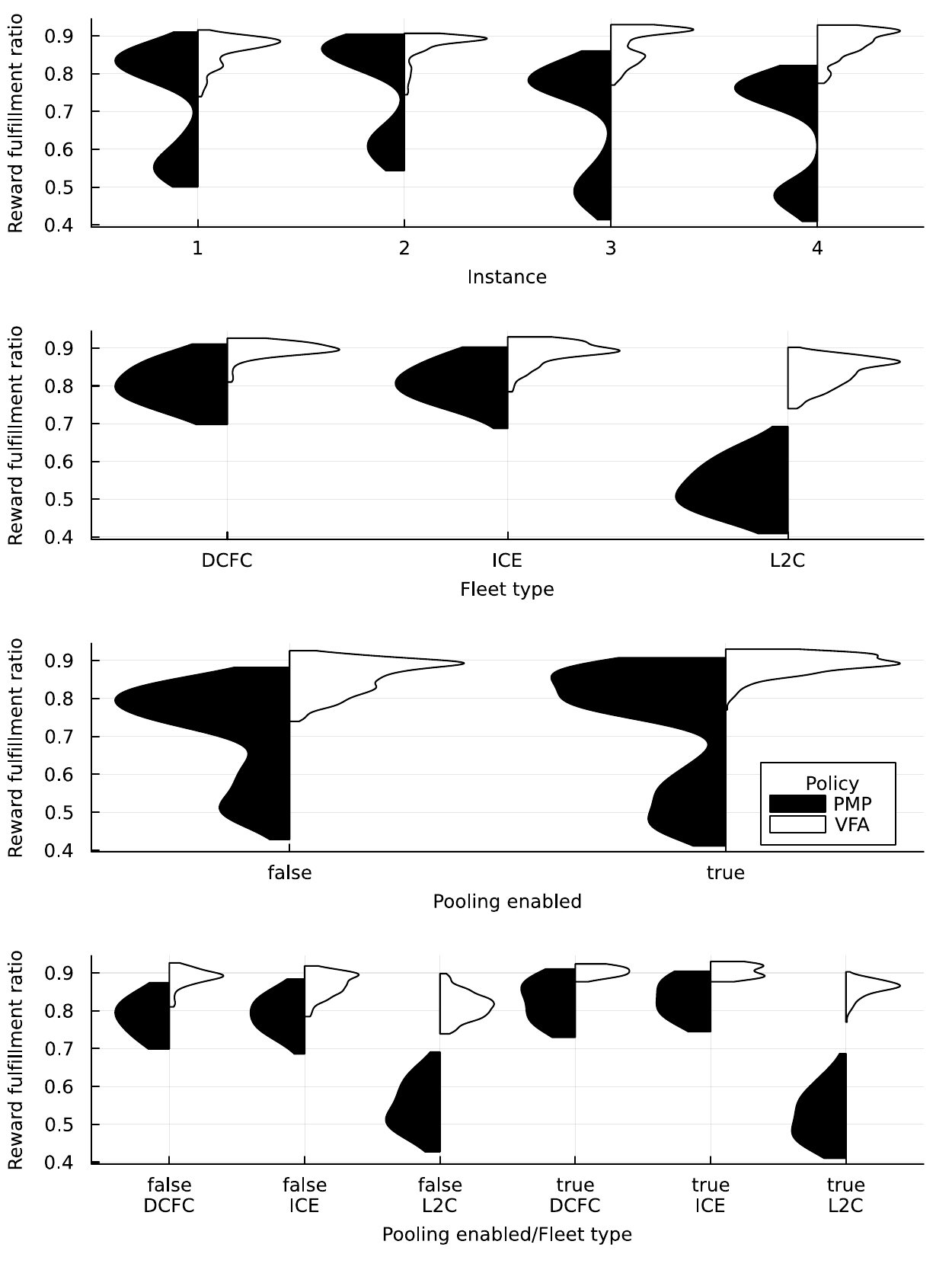}
  \caption{Kernel density plots of the \gls{rfr} obtained by the \gls{pm} and \gls{vfa} policies across different instances and problem settings.}
  \label{fig:rfr_violins}
\end{figure}

We report the \gls{rfr}, which we define as the ratio between the reward obtained by a policy and the sum of the rewards of all requests in the sample path. \Cref{fig:rfr_violins} shows kernel density plots of the \gls{rfr} obtained by the \gls{pm} and \gls{vfa} policies across different instances and problem settings. \Cref{sec:tables} contains tables of reward and \gls{rfr} statistics for each instance, problem setting, and policy.

\newcommand{\result}[3]{#2 #3}
\newcommand{\mediqr}{median \gls{rfr} (\gls{iqr})}

\paragraph{Result 1: The \gls{vfa} policy performs significantly better than the \gls{pm} policy.}
\result{The \gls{vfa} policy performs significantly better than the \gls{pm} policy.}{The \mediqr{} of the \gls{pm} and \gls{vfa} policies is \PMRfrMedianIqr{} and \VFARfrMedianIqr{}, respectively. \Cref{fig:rfr_violins} shows that the difference in performance cannot be fully explained by controlling for instance or problem characteristics. Nevertheless, the difference is more pronounced in instances 3 and 4, and when using \gls{l2c} \glspl{EV}.}{The advantage of the \gls{vfa} policy is consistent with the literature. Since the \gls{pm} policy is clearly suboptimal, we ignore it for the remainder of this analysis.}

\paragraph{Result 2: Both \gls{ICE} vehicles and \gls{DCFC} \glspl{EV} achieve the highest \gls{rfr}, while \gls{l2c} \glspl{EV} perform worse.}
\result{Both \gls{ICE} and \gls{DCFC} vehicles achieve the highest \gls{rfr}, while \gls{l2c} perform worse.}{The \mediqr{} of \gls{ICE} and \gls{DCFC} vehicles is \VFARfrMedianIqrICE{} and \VFARfrMedianIqrDCFC{}, respectively. By contrast, the \mediqr{} of \gls{l2c} vehicles is \VFARfrMedianIqrLTwoC. This is illustrated in the fleet types subplot in \Cref{fig:rfr_violins}.}{
The relatively weak performance of \gls{l2c} vehicles is expected, because one would expect the vehicle range and the recharging or refueling rate $\symbolRechargingRate$ to have a positive correlation with the \gls{rfr}. 
Interestingly, the performance of \gls{DCFC} vehicles is on par with \gls{ICE} vehicles despite the significant differences in their range and $\symbolRechargingRate$ values.
One possible explanation for this is that the comparatively small distances of metropolitan areas such as \gls{NYC} remove the significance of the vehicle range.}

\paragraph{Result 3: Pooling enables the \gls{rhs} to achieve high rewards more consistently.}
\result{Pooling enables the \gls{rhs} to achieve high rewards more consistently.
}{
When pooling is not enabled, the median \gls{rfr} is \VFARfrMedianPoolingfalse{}, which improves to \VFARfrMedianPoolingtrue{} with pooling. More noteworthy than the improvement in the median \gls{rfr} is the change in the shape of the \gls{rfr} distribution, as visualized in the third subplot from the top in \Cref{fig:rfr_violins}. The change in the shape of the distribution is associated with a decrease in the \gls{iqr} from \VFARfrIqrPoolingfalse{} to \VFARfrIqrPoolingtrue{}. When we focus on the \gls{DCFC} and \gls{ICE} settings (see Result 4 below for a detailed analysis of the \gls{l2c} setting), the bottom subplot shows that pooling makes their \gls{rfr} distributions especially consistent, without significantly affecting their best-case \gls{rfr}.}{One possible explanation for this is the following. Without pooling, the \gls{vfa} policy already achieves the maximum \gls{rfr} in some ``easy'' sample paths, and this maximum is around 90\% \gls{rfr}. The remaining requests in the gap to 100\% \gls{rfr} are not profitable because their reward is low, they take very long, or they connect low-demand areas. As a result, enabling pooling has no effect on an ``easy'' sample path. However, there are other, more difficult sample paths where the policy needs pooling to achieve maximum \gls{rfr}.
}

\begin{figure}[H]
    \begin{subfigure}{0.49\textwidth}
        \includegraphics[width=0.975\textwidth]{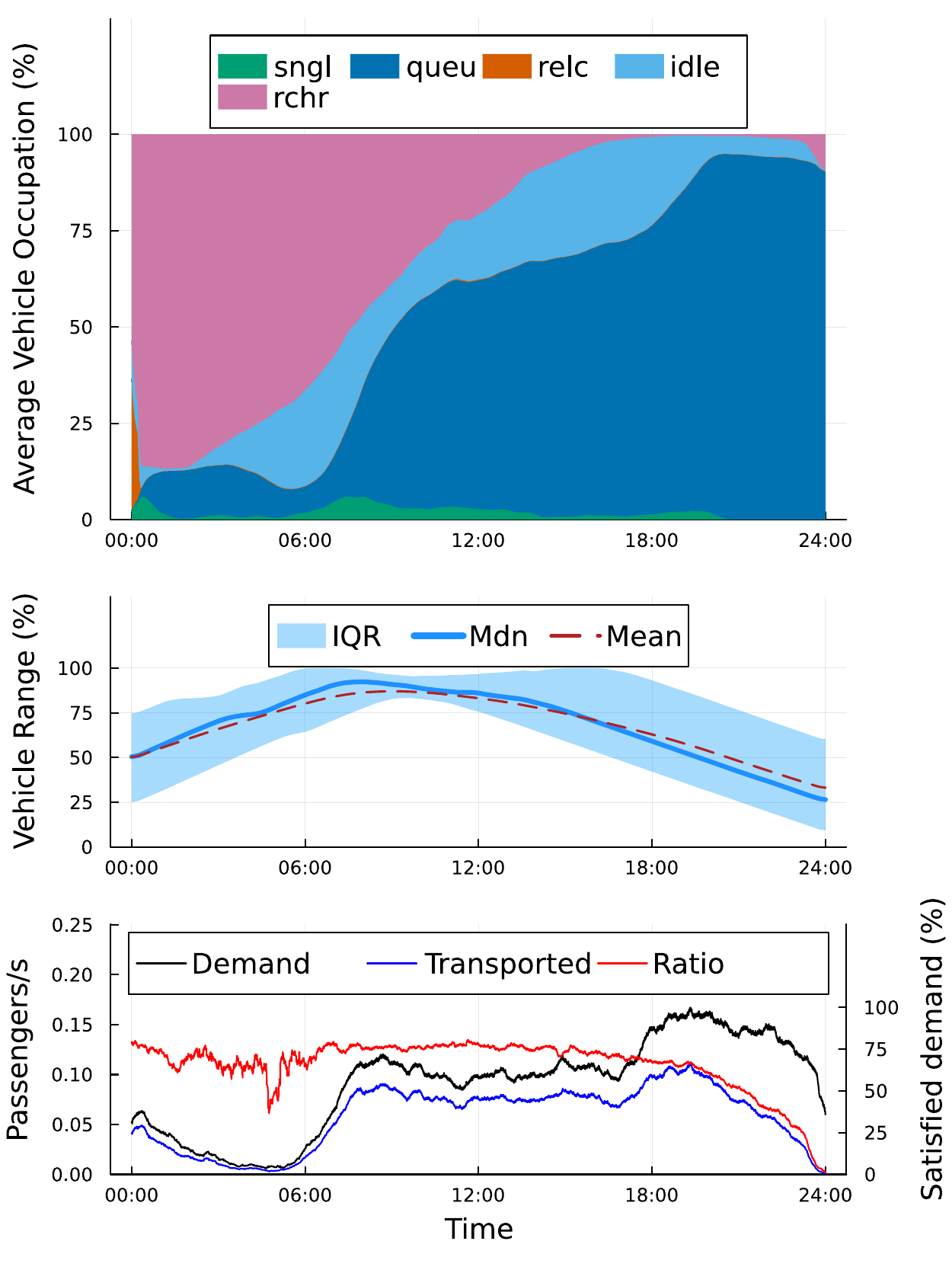}
    \end{subfigure}
    \hfill
    \begin{subfigure}{0.49\textwidth}
        \includegraphics[width=0.975\textwidth]{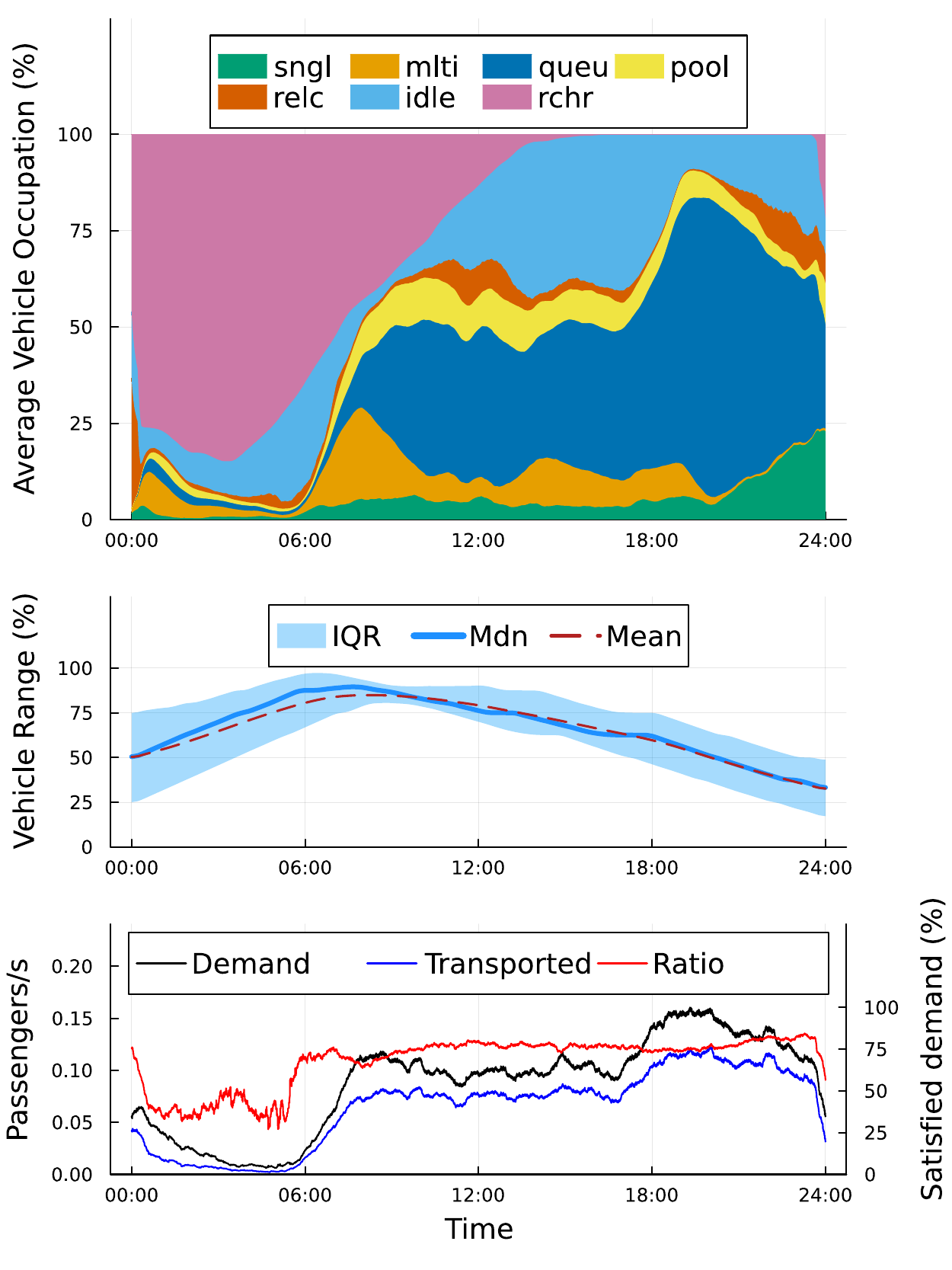}
    \end{subfigure}
    \caption{Behavior of the \gls{vfa} policy in instance 1, using  \gls{l2c} vehicles, without (left) and with pooling (right), averaged across 30 sample paths.}
    \label{fig:ORD-4-vfa-l2c}
\end{figure}

\paragraph{Result 4: Pooling alleviates the disadvantage of \gls{l2c} vehicles.}
\result{Pooling alleviates the disadvantage of \gls{l2c} vehicles.}{Pooling improves the \mediqr{} of \gls{l2c} fleets from  \VFARfrMedianIqrPoolingfalseLTwoC{} to \VFARfrMedianIqrPoolingtrueLTwoC. By contrast, the \mediqr{} of \gls{DCFC} vehicles merely improves from \VFARfrMedianIqrPoolingfalseDCFC{} to \VFARfrMedianIqrPoolingtrueDCFC, with similar numbers for \gls{ICE} vehicles. \Cref{fig:ORD-4-vfa-l2c} illustrates the influence of pooling on \gls{l2c} fleets. Comparing the top subplots shows that the pooling (right) policy executes multi-trip and pool decisions, which the non-pooling (left) policy does not. However, this is not the only difference: The pooling policy also executes significantly more relocate decisions. Comparing the middle subplots shows that towards the end of the day, the pooling (right) fleet's average range is higher than the non-pooling (left) fleet's average range, even though both fleets have similar recharging behavior. The bottom subplots show that the pooling (left) fleet is able to transport much more passengers in the high-demand hours between 6 p.m.\ and midnight than the non-pooling (right) fleet.}{The data suggests that the pooling decisions allow the \gls{rhs} to serve the same passengers with fewer vehicles and relocate the others to more promising locations, while using the range of all vehicles more efficiently. The effect of this on the total reward is most noticeable when the need to make efficient use of the vehicles' range is highest. This makes the \gls{l2c} setting stand out.}

\section{Conclusion}\label{sec:conclusion}
Optimal control of ride-hailing systems is a challenging problem that involves several interdependent decisions, which must be made in a dynamic and stochastic environment while adhering to various operational constraints. Although the problem has been studied extensively in the literature, even in comparable settings, studies differ significantly
in the types of problem features, decisions, constraints, and performance metrics they consider. This lack of consistency makes it difficult to draw generalizable conclusions and undermines the practical relevance of the findings. This is exacerbated by the fact that many studies do not publish their benchmark instances, making it hard to replicate or build upon their findings. 

To address this gap, we proposed a modular, unified modeling framework that integrates various problem features, including pooling, repositioning, \glspl{EV} and \glspl{ICE} fleets, as well as various constraints regarding customer response time, waiting time, and detour time. To solve the problem, we designed a \gls{pm} policy that charges (or fuels) vehicles when their batteries (or fuel tanks) are below a certain threshold, and a \gls{vfa} policy that captures the long-term impact of current decisions. Both policies utilize an efficient procedure for enumerating all feasible vehicle-to-request assignments. Furthermore, we introduced a scalable technique that introduces randomness in the decision-making process, as well as monotonicity properties of the value function, to effectively address the exploration-exploitation tradeoff in \gls{adp}.
We also created reusable benchmark instances that capture a range of spatial structures and demand distributions, based on real-world data from \gls{NYC}.
Our computational results provide insights into previously unstudied interactions between problem characteristics. Specifically, we found no significant difference between revenue generated by \gls{ICE} fleets and fast-charging \gls{EV} fleets, but both significantly outperformed slow-charging \gls{EV} fleets. We also found that pooling increases the revenue, and reduces revenue variability, for all fleet types.

Our unified modeling framework can serve as a foundation for several promising future research directions. Relevant extensions include the integration of pre-scheduled trips (as seen in airport shuttles or paratransit services),
heterogeneous fleets (with varying seat capacities or energy sources), stochastic travel times and congestion, uncertain passenger headcounts and destinations that are only revealed upon pickup, infrastructure limitations (e.g., scarce parking at key locations and queuing at charging or fueling stations), and explicit modeling of customers' varying preferences (e.g., willingness to pool and tolerance for detours).
On the solution side, future work may benefit from spatial aggregation techniques that enhances scalability without compromising solution quality.
Additionally, exploring alternative policy architectures, such as regression-based value function approximations, direct lookahead policies, and more expressive representations using deep learning models, including large-scale transformers, all present promising avenues for further investigation.

\bibliographystyle{apalike}
\bibliography{00_bibliography}

\clearpage
\begin{appendix}
This appendix is structured as follows.
\Cref{sec:transition_function_details} defines the transition function of our model.
\Cref{sec:policies_details} explains how to derive and implement our \gls{vfa} policy, and establishes some definitions. \Cref{sec:proof} gives a proof of correctness of the multi-trip decision enumeration algorithm that we use in our \gls{vfa} policy.
\Cref{sec:experiments_details} contains implementation details and additional results of our numerical studies.

\section{Definition of the transition function}
\label{sec:transition_function_details}

To define $a^{M}(a,d)$, we first introduce the following notation:
\begin{itemize}
    \item $\tau(v,w)$ is the time it takes to traverse the shortest path from node $u\in\V$ to node $v\in\V$.
    \item $\tau(p)$ is the time it takes to traverse the entire path $p\in P$.
    \item $\sigma(u,v)$ is the driving range required to travel the shortest path from node $u\in\V$ to node $v\in\V$.
    \item $\sigma(p)$ is the driving range required to traverse path $p\in P$.
    \item $\pathend$ is the last node in the path $p\in P$.
    \item $\ell(l_a)$ is the time it takes a vehicle $a$ to recharge its battery level from $l_a$ to $l^{\max}$.
    \item $f_a$ is the first node on the shortest path from $o_a$ to $d_a$ that vehicle $a\in\A$ reaches after $\lfloor\taua\rfloor + 1$.
    \item $m_a(h) = \max\left\{\taua + h, \lfloor\taua\rfloor + 1\right\}$, where $h$ is a time span.
\end{itemize}
The function $a^{M}(a,d)$ is defined for every attribute $a \in \A$, decision $d \in \D(a)$, and decision epoch $\tInSetTimeSteps$. Below we specify the definition of $a^{M}(a,d)$ for every type of decision $d \in \D(a)$.

\paragraph{Continue:} A vehicle $a\in\Aoccu$ that gets assigned to a decision $\deccont \in \Deccont_a$ at time $t$ spends the time between its actionable time $t_a$ and the next decision epoch $t+1$ driving on the shortest path from its current location $o_a$ to its destination $d_a$. If $t_a \geq t+1$, the future attribute is equal to $a$. If $t_a > t+1$ and the vehicle reaches $d_a$ before time $\lfloor\taua\rfloor + 1$, it parks until then. Otherwise, it becomes actionable in the first location between $o_b$ and $d_b$ that it reaches after $\lfloor\taua\rfloor +1$. The future attribute is
\begin{itemize}
    \item $a^{M}(a,\deccont) = a$ if $\taua \geq t+1$,
    \item $a^{M}(a,\deccont) = \left(d_a, \; d_a, \; l_a - \sigma(o_a, d_a), \; n_{\max}, \; \lfloor\taua\rfloor + 1\right)^{\intercal}$ if $\taua < t+1 \wedge \taua+\tau(o_a, d_a) \leq \lfloor\taua\rfloor + 1$, and 
    \item $a^{M}(a,\deccont) = \left(f_a, \; d_a, \; l_a - \sigma(o_a, f_a), \; n_a, \; \taua+\tau(o_a, f_a)\right)^{\intercal}$ otherwise.
\end{itemize}

\paragraph{Single trip:}
A vehicle $a \in \Aempt$ that gets assigned to a decision $\decs_b \in \Decs_a$ to serve a request $b \in \B$ starts moving from its current location $o_a$ to the origin $o_b$ to pickup the passengers. After pickup, the vehicle acts like under a \emph{continue} decision. Its future attribute is $a^{M}(a,\decs_b) = a^{M}(a',\deccont)$, where
\begin{equation*}
     a' = \begin{pmatrix}o_b \\ d_b \\ l_a - \sigma(o_a, o_b) \\ n^{\max} - n_b \\ t_a + \tau(o_b, d_b)\end{pmatrix}
\end{equation*}

\paragraph{Multi-trip and pool:}
A vehicle $a \in \A$ that gets assigned to a decision $\decm_{Bp} \in \Decm_a$ or $\decp_{bp} \in \Decp_a$ must drive the entire path $p$ to serve requests $B$, or $b$, respectively. If this does not take until the next decision epoch, the vehicle idles at the last drop off location until then. The future attribute of these vehicles is $a^{M}(a,\decm_{Bp}) = a^{M}(a,\decp_{bp}) = \left(\pathend, \pathend, l_a - \sigma(p), n^{\max}, m_a(\tau(p)) \right)^\intercal$.
\paragraph{Queue:}
A vehicle $a \in \Aoccu$ that gets assigned to a decision $\decq_b \in \Decq_a$ drives to $d_a$, drops off its passengers, then drives directly to the origin $o_b$ of request $b \in \B$ and picks up the new passengers. After pickup, the vehicle acts like under a \emph{continue} decision. Its future attribute is $a^{M}(a,\decq_b) = a^{M}(a',\deccont)$, where
\begin{equation*}
     a' = \begin{pmatrix}o_b \\ d_b \\ l_a - \sigma(o_a, d_a) - \sigma(d_a, o_b) \\ n^{\max} - n_b \\ t_a + \tau(o_a, d_a) + \tau(d_a, o_b)\end{pmatrix}
\end{equation*}

\paragraph{Relocate and recharge:} 
The future attribute of a vehicle $a\in\Aempt$ that gets assigned to a relocate or recharge decision is $a^{M}(a,\decrel_v) = \left(v, v, l_a - \sigma(o_a,v), n^{\max}, t+1\right)^{\intercal}$
and
$
a^{M}(a,\decrel_v) = \left(o_a, o_a, l^{\max}, n^{\max}, m_a(\ell(l_a)) \right)^{\intercal},
$
respectively.

\paragraph{Idle:}
A vehicle $a\in\Aempt$ that gets assigned to an idle decision parks in its current location until the next decision epoch. If the actionable time $t_a \geq t+1$, the future attribute stays the same. This allows $t$ to catch up to $t_a$. The future attribute of such a vehicle is $a^{M}(a,\decidle) = a$ if $\taua \geq t+1$, and 
$a^{M}(a,\decidle) =
\left(o_a, \; d_a, \; l_a, \; n_{\max}, \; \lfloor\taua\rfloor + 1\right)^{\intercal}
$
otherwise.

\section{Policies Supplement}
\label{sec:policies_details}

This section is structured as follows.
\Cref{sec:derivation} shows the derivation of the \cite{simao:2009} linear \gls{vfa} in the context of our model.
\Cref{sec:forward_adp_algorithm} specifies the details of the Forward \gls{adp} algorithm. 
\Cref{sec:monotonicity} discusses the monotonicity properties we use to learn our approximate value function.

\subsection{Derivation of linearized value function}
\label{sec:derivation}
The value function is defined recursively as
\begin{equation}\label{eqx:bellman_equation}
 V(S_t) = \max_{x_t \in \X(S_t)} \big\{ C(x_t) + \mathbb{E}[V(S_{t+1}) \mid S_t, x_t] \big\},
\end{equation}
for $\tInSetTimeSteps$ and $V(S_{T+1})=0$.

The first difficulty in dealing with the Bellman equation~\eqref{eqx:bellman_equation} is the expectation term, which is taken with respect to the random exogenous demand $\hat D_{t+1}$. We can move the expectation outside the max operator by using the post-decision state $S_t^x$ and replacing $\mathop{\mathbb{E}} [V(S_{t+1}) \; | \; S_t, x_t]$ with $V^x(S_t^x)$ in \eqref{eqx:bellman_equation}, where $V^x(S_t^x) = \mathop{\mathbb{E}} [V(S_{t+1})|S_t^x]$.
Readers familiar with other reinforcement learning algorithms may find it helpful to note that approximating $V^x(S_t^x)$ is essentially the same as Q-learning, except that we estimate the value of post-decision states rather than state-action pairs.

The second difficulty in dealing with the Bellman equation~\eqref{eqx:bellman_equation} is the exponentially large number of possible states, which makes it intractable to obtain the exact value of $V^x(S_t^x)$. A typical approach to deal with this difficulty is to approximate $V^x(S_t^x)$ by a function $\bar{V}^x(S^{x}_t)$ that is easier to compute. We use a linear approximation \citep{simao:2009} given by
\begin{equation}\label{eqx:linear_approximation}
 \bar{V}^x(S^{x}_t) = \sum_{a \in \A} \bar{v}^{R}_{a}R^{x}_{ta} + \sum_{b\in\B} \bar{v}^{D}_{tb}D^{x}_{tb}
\end{equation}
where
\begin{itemize}
 \item $\bar{v}^{R}_{a}= \frac{\partial V(S_t)}{\partial R_{ta}}$ is the marginal value of a resource with attribute $a\in \A$ at any time $\tInSetTimeSteps$, and
 \item $\bar{v}^{D}_{tb}= \frac{\partial V(S_t)}{\partial D_{tb}}$ the marginal value of a demand with attribute $b\in \B$ at time $t$.
\end{itemize}
Note that we omit the time index $t$ in $\bar{v}^{R}_{a}$ since the vehicle's actionable time is already encoded in the attribute $a$, and we estimate a single marginal value $\bar{v}^{R}_{a}$ across all $\tInSetTimeSteps$.
By substituting $\bar{V}^x(S^{x}_t)$ into \Cref{eqx:bellman_equation}, we obtain
\begin{align*}
 \hat{V}(S_t) \approx
 & \max_{x_t \in \X(S_t)}\left\{C(x_t) + \sum_{a \in \A} \bar{v}^{R}_{a}R^{x}_{ta} + \sum_{b\in\B} \bar{v}^{D}_{tb}D^{x}_{tb} \right\}\\
 =
 & \max_{x_t \in \X(S_t)}\left\{ \sum_{a \in \A}\sum_{d \in \D(a)} c_{tad} x_{tad}+ \sum_{a' \in \A} \bar{v}^{R}_{a'}\sum_{a \in \A}\sum_{d \in \D(a)} x_{tad} \cdot \mathbf{1}_{\{a^{M}(a,d)=a'\}}+\sum_{b\in\B} \bar{v}^{D}_{tb}D^{x}_{tb}\right\}\\
 =
 & \max_{x_t \in \X(S_t)}\left\{\sum_{a \in \A}\sum_{d \in \D(a)} \left(c_{tad}+ \sum_{a' \in \A} \bar{v}^{R}_{a'}\cdot \mathbf{1}_{\{a^{M}(a,d)=a'\}} \right)x_{tad}+\sum_{b\in\B} \bar{v}^{D}_{tb}D^{x}_{tb}\right\}\\
 =
 & \max_{x_t \in \X(S_t)}\left\{ \sum_{a \in \A}\sum_{d \in \D(a)} \left(c_{tad} + \bar{v}^{R}_{a^{M}(a,d)} \right)x_{tad}+\sum_{b\in\B} \bar{v}^{D}_{tb}D^{x}_{tb}\right\},
\end{align*}
which gives us our \gls{vfa} policy
\begin{align}\label{eqx:vfa_policy}
\pi^\text{VFA}(S_t) = \argmax_{x_t \in \X(S_t)}\left\{ \sum_{a \in \A}\sum_{d \in \D(a)} \left(c_{tad} + \bar{v}^{R}_{a^{M}(a,d)} \right)x_{tad}+\sum_{b\in\B} \bar{v}^{D}_{tb}D^{x}_{tb}\right\}.
\end{align}

\subsection{Forward \gls{adp} algorithm}
\label{sec:forward_adp_algorithm}
\begin{algorithm}[htbp]
    \caption{Forward Approximate Dynamic Programming}\label{alg:forward_adp}
    \begin{algorithmic}[1]
        \REQUIRE
        Initial resource values $\bar{v}^{R,0}_{a}$ for all $a \in \A$.
        Initial demand values $\bar{v}^{D,0}_{tb}$ for all $b \in \B, \tInSetTimeSteps$.
        Initial states $S^n_0$, sample paths $(\hat D^n_1, \hat D^n_2, \dots, \hat D^n_T),$ and learning rates $\alpha_n$  for all $n \in 1, \dots, N$.
        \FOR{$n\in\{1,\dots, N\}$}
        \FOR{$\tInSetTimeSteps$}
        \STATE \label{ln:adp_iteration} Solve the LP \[\max_{x^n_t \in \X^\mathrm{LP}(S^n_t)}\left\{ C(x^n_t) + \sum_{a \in \A} \bar{v}^{R,n-1}_{a}R^{x}_{ta}+\sum_{b\in\B} \bar{v}^{D,n-1}_{tb}D^{x}_{tb}\right\}\] to obtain $x^n_t$, $\hat{v}^{R,n}_{ta}$ for all $a\in\A$ and $\hat{v}^{D,n}_{tb}$ for all $b\in\B$.
        \STATE \label{ln:update} Update the lookup tables \begin{align*}
            \bar{v}^{R,n}_{a} &:= (1-\alpha_n) \bar{v}^{R,n-1}_{a} + \alpha_n \hat{v}^{R,n}_{ta},\\
            \bar{v}^{D,n}_{tb} &:=(1-\alpha_n) \bar{v}^{D,n-1}_{tb} + \alpha_n \hat{v}^{D,n}_{tb}.
        \end{align*}
        \STATE Transition to the next state $S^n_{t+1} = S^{M}(S^n_t, x^n_t, \hat D^n_{t+1})$.
        \ENDFOR
        \ENDFOR
        \RETURN $\bar{v}^{R,N}_{a}$ and $\bar{v}^{D,N}_{tb}$ for $\tInSetTimeSteps$.
    \end{algorithmic}
\end{algorithm}

We adapt the well-known Forward \gls{adp} algorithm (see \Cref{alg:forward_adp}) to calculate estimates of $\bar{v}^{R}_{a}$ for every $a \in \A$ and $\bar{v}^{D}_{tb}$ for every $b\in \B$ and $\tInSetTimeSteps$, which we use in our \gls{vfa} policy \eqref{eqx:vfa_policy}.
The idea of the Forward \gls{adp} algorithm is to start from an initial approximation $\bar{v}^{R,0}_{a}$ and $\bar{v}^{D,0}_{tb}$, simulate $N$ sample paths of the system, and iteratively improve these approximations by solving an \gls{lp} at every time step of every sample path.
Specifically, in each iteration $n \in \{1, \dots, N\}$, we simulate a sample path of the system from time $t=0$ to $t=T$. Then, at every time step $\tInSetTimeSteps$, we evaluate a surrogate policy by removing all multi-trip decisions from the feasible region $\X(S_t)$, relaxing the integrality constraint, and solving an \gls{lp} that uses the current approximations $\bar{v}^{R,n-1}_{a}$ and $\bar{v}^{D,n-1}_{tb}$ (Line \ref{ln:adp_iteration} in \Cref{alg:forward_adp}). Solving this \gls{lp} is guaranteed to yield an integer solution, due to the total unimodularity of the constraint matrix and the integrality of the right-hand side. Therefore, solving this \gls{lp} we can obtain the shadow prices
\begin{itemize}
    \item $\hat{v}^{R,n}_{ta} = \frac{\partial \hat V(S^n_t)}{\partial R_{ta}}$ for the constraint of vehicle attribute $a\in\A$, and 
    \item $\hat{v}^{D,n}_{tb} = \frac{\partial \hat V(S^n_t)}{\partial R_{ta}}$ for the constraint of request attribute $b\in\B$. 
\end{itemize}
where $S_{t}^n$ is the system state at time $t$ in the $n$-th iteration, and $\hat{v}^{R,n}_{ta}$ and $\hat{v}^{D,n}_{tb}$ are estimates of the slope of $\hat V(S^n_t)$ with respect to $R_{ta}$ and $D_{tb}$ at the point $R^{x}_{ta}$ and $D^{x}_{tb}$, respectively. Note that we retain the time index $t$ in $\hat{v}^{R,n}_{ta}$ to emphasize that the partial derivative is taken with respect to the resource vector $R_{ta}$ at time $t$. 
Finally, we use these estimates to update the current approximations $\bar{v}^{R,n-1}_{a}$ and $\bar{v}^{D,n-1}_{tb}$ as shown in Line \ref{ln:update} of \Cref{alg:forward_adp}, transition to the next state $S^n_{t+1}$, and repeat this process until we reach the end of the sample path.

\subsection{Monotonicity of the value function}
\label{sec:monotonicity}

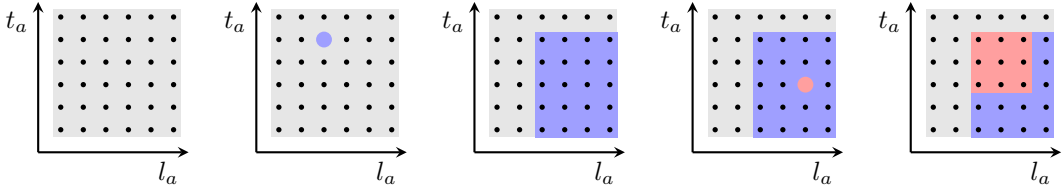
\begin{figure}
    \centering
    \begin{tikzpicture}[scale=1, transform shape,>=stealth, font=\sffamily]

\definecolor{gatorgray}{gray}{0.9}     
\colorlet{gatororange}{red!50}           
\colorlet{gatorblue}{blue!50}
\newcommand{\drawdots}{%
\foreach \i in {0.3,0.6,...,2.1} {
    \foreach \j in {0.3,0.6,...,2.1} {
        \fill[black] (\i,\j) circle (1pt);
    }
}
}
\newcommand{\drawaxes}{
\draw[->,thick] (0,0) -- (2,0) node[below left] {$l_a$};
\draw[->,thick] (0,0) -- (0,2) node[below left] {$t_a$};
}

\newcommand{\drawgray}{%
\drawaxes
\fill[gatorgray] (0.2,0.2) rectangle (1.9,1.9);
}
\newcommand{\drawblue}{%
\fill[gatorblue!75] (0.8,1.6) rectangle (1.9,0.2);
}
\begin{scope}[xshift=0cm, yshift=0cm]
\drawgray
\drawdots
\end{scope}

\begin{scope}[xshift=2.9cm, yshift=0cm]
\drawgray
\drawdots
\fill[gatorblue!75] (0.9,1.5) circle (3pt);
\end{scope}
\begin{scope}[xshift=5.8cm, yshift=0cm]
\drawgray
\drawblue
\drawdots
\end{scope}
\begin{scope}[xshift=8.7cm, yshift=0cm]
\drawgray
\drawblue
\drawdots
\fill[gatororange!75] (1.5,0.9) circle (3pt);
\end{scope}

\begin{scope}[xshift=11.6cm, yshift=0cm]
\drawgray
\drawblue
\fill[gatororange!75] (0.8,1.6) rectangle (1.6,0.8);
\drawdots
\end{scope}
\end{tikzpicture}
    \caption{Monotonicity properties are used to update a value function twice. The $o_a$, $d_a$ and $n_a$ dimensions are not drawn to allow for a two-dimensional representation. Initially, the value function is at a low value (left). The first update increases it to a high value for all attributes $a$ with $a \preccurlyeq a'$ (center). The second update decreases it to a medium value for all attributes $a$ with $a'' \preccurlyeq a$ (right). Figure adapted from \citet{powellReinforcementLearningStochastic2022}.}
    \label{fig:monotonicity}
\end{figure}

This section extends Section 4.3.3 of the paper. We first repeat the relevant definition:

\partialOrderDefinition{x}

Definition \ref{defx:partial_order} allows us to show three monotonicity properties with respect to an optimal value function for vehicle attributes, $v^{R*}_a$.
The first property applies when we compare two vehicles $a\preccurlyeq a'$ and $a'$ has a higher range than $a$. In this case, the higher-range vehicle cannot have a lower value than the lower-range vehicle. This is because the higher-range vehicle can serve the same trips as the lower-range vehicle and therefore the higher-range vehicle can obtain the same reward as the lower-range vehicle.
The second property applies when we compare two vehicles $a\preccurlyeq a'$ and $a'$ has a higher capacity than $a$. In this case, the higher-capacity vehicle cannot have a lower value because it can serve the same trips as the lower-capacity vehicle.
The third property applies when we compare two vehicles $a\preccurlyeq a'$ and $a'$ can start performing new tasks (i.e., is actionable) at an earlier time than $a$. In this case, the earlier vehicle cannot have a lower value because it has more time to serve travel requests and therefore can earn the same reward as the later vehicle.

\monotonicityProposition{x}

\begin{proof}
    For every pair of vehicle attributes $a,a' \in \A$, every decision $d\in \D(a) \cap \D(a')$, and every time step $\tInSetTimeSteps$, the statement $
    a\preccurlyeq a' \wedge a^M(a,d) \preccurlyeq a^M(a',d) \wedge c_{tad} \leq c_{ta'd}$ is satisfied by definition of $a^M$ and $c_{tad}$.
    All three monotonic properties can be proved using the same backward induction argument as in~\cite{alkanj:2020}. 
\end{proof}
We note that the following two assumptions made by \citet{alkanj:2020} for their monotonicity properties are not necessary for us since the resource vector $R_t$ is deterministic in our model. The first assumption is that for every $a, a' \in \A$, if $a \preccurlyeq a'$ is satisfied in some $S^x_t$, it is also satisfied in $S_{t+1}$. The second assumption is that the demand component and the vehicle component of the exogenous information are stochastically independent.

We leverage the monotonicity of ${v}^{R*}_a$ to address the exploration-exploitation trade-off, and to accelerate learning $\bar v^{R,n}_a$, as follows.
We initialize the lookup table such that it is monotonic in all dimensions (e.g. $\bar v^{R,0}_a=0,\, \forall a\in\A$).
After updating the lookup table (lines 7 and 8 of Algorithm \ref{alg:forward_adp}), $\bar v^{R,n}_a$ may violate some monotonicity property.
To restore monotonicity, we force all entries in the lookup table that violate monotonicity to a larger or smaller value so that monotonicity is restored. See Figure \ref{fig:monotonicity} for an illustration.

\FloatBarrier
\section{Proof of multi-trip decision enumeration algorithm}
\label{sec:proof}

We first repeat the relevant definitions, then state the proposition, and finally give the proof.
\detourPenaltyMonoDefinition{x}
\enumerationDefinitions{x}
\enumerationProposition{x}
\begin{proof}
    Suppose, for the sake of contradiction, that the proposition is false. 
    The negation of the proposition is as follows.
    
    There exists a time step $\tInSetTimeSteps$ and a state $S_t$ where a vehicle attribute $a\in\Aempt$ satisfies $R_{ta} > 0$ and where the nonempty set of requests $B\subseteq\mathcal{B}$ satisfy $\forall b\in B.D_{tb} > 0$ and $|P_{aB}| > 0$ such that
    (i) the approximate value function used by policy $\pi^\mathrm{VFA}$ satisfies \Cref{prop:monotonicity},
    (ii) $\detourPenalty{}$ is monotonically increasing in the path length, and 
    (iii) there exists an \gls{ldsps} for $(a,B)$ that is not a-B-$\pi^\mathrm{VFA}$-sufficient.
    
    Statement (iii) is equivalent to the following statement:
    There exists a set of paths $P_{aB}^* \subseteq P_{aB}$ such that
    (iv) the following statement is true for all requests $b\in B$: Either $P_{aB}^*$ contains exactly one of the shortest paths in $P_{aB}$ that end in the destination of $b$, $d_b$, or no path in $P_{aB}$ ends in $d_b$, and
    (v) for all solutions $x_t$ in the set of optimal solutions $\pi^\mathrm{VFA}(S_t)$, there exists a decision variable $x_{tad}$ with $d\in \{\decm_{Bp} \in \Decm_a \mid p\in P_{aB}\setminus P^*_{aB}\}$ that is not equal to zero.
    
    Analyzing Statement (v), let $p \in P_{aB}\setminus P^*_{aB}$ be the path associated with $x_{tad}$, i.e., $d = \decm_{Bp}$. Since $B$ is nonempty, there exists a request $b\in B$ that is dropped off last in $p$. Using Statement (iv), this implies that the set $P_{aB}^*$ is nonempty, so there exists a shortest path $p^*\in P_{aB}^*$ with $\tau(p^*) \leq \tau(p)$.
    By \Cref{defx:partial_order}, this is equivalent to $a^M(a,d) \preccurlyeq a^M(a,d^*)$, where $d^* = \decm_{Bp^*}$. Because of Statements (i) and (ii), this implies that 
     \begin{equation}\label{eqx:dsmallerdprime}
        \bar{v}^{R}_{a^M(a,d)} \leq \bar{v}^{R}_{a^M(a,d^*)} \text{\quad{}and\quad{}} \detourPenalty{}(a,d) \leq \detourPenalty{}(a,d^*).
    \end{equation}
    
    We have shown that the supposition implies that both Statement (v) and  \eqref{eqx:dsmallerdprime} are true. We will now show that the supposition also implies that either Statement (v) or \eqref{eqx:dsmallerdprime} is false.
    
    By construction, the optimization problem \eqref{eqx:vfa_policy} contains the decision variable $x_{tad^*}$, and the objective function contains both
    \begin{equation*}
        \Bigl (c_{tad} + \bar{v}^{R}_{a^M(a,d)} \Bigr )x_{tad}\text{\quad{}and\quad{}}\Bigl ( c_{tad^*} + \bar{v}^{R}_{a^M(a,d^*)} \Bigr ) x_{tad^*}
    \end{equation*}as summands.
    Because of the nonnegativity constraint in $\X$, $x_{tad} \geq 0$ and $x_{tad^*} \geq 0$ in all solutions $\pi^\mathrm{VFA}(S_t)$. 
    Using Statement (v), this implies that $x_{tad} > 0$ in all solutions.
    Two cases remain: $x_{tad^*} > 0$ or $x_{tad^*} = 0$.
    
    If $x_{tad^*} > 0$, the objective coefficients of both variables must be equal. This implies that there is an alternative optimal solution with $x_{tad} = 0$. This contradicts Statement (v).
    
    If $x_{tad^*} = 0$, this is equivalent to 
     \begin{align*}
        c_{tad} + \bar{v}^{R}_{a^M(a,d)}
        &> c_{tad^*} + \bar{v}^{R}_{a^M(a,d^*)}\\
        \Leftrightarrow \quad
        \detourPenalty{}(a,d) + \bar{v}^{R}_{a^M(a,d)}
        &>\detourPenalty{}(a,d^*) + \bar{v}^{R}_{a^M(a,d^*)}
    \end{align*}
    which contradicts \eqref{eqx:dsmallerdprime}.
\end{proof}

\section{Supplement to numerical experiments}
\label{sec:experiments_details}
This section outlines the process of creating the benchmark instances we use in our numerical experiments, discusses the problem settings parameters, explains relevant implementation details, and gives detailed results.

\subsection{Street networks} 
\label{sec:street_networks}

The 2009 \gls{TLC} data provides the latitude and longitude coordinates of every trip's origin and destination. This allows us to adopt a fine-grained representation of the street network where we discretize \gls{NYC} into rectangles of equal size, each approximately \qty{215}{\meter} wide by \qty{280}{\meter} high. We define a directed graph that has one vertex per rectangle and that contains an arc $(u, v)$ if and only if vertices $u$ and $v$ represent rectangles that are directly connected by a street according to OpenStreetMap data. 
The \gls{TLC} data does not contain a trip between each adjacent pair of vertices that we have in our graph. To define the driving duration of each arc, we fit the simple linear regression model $t = \beta s$, where $t$ is the trip duration in seconds and $s$ is the haversine distance in meters between the origin and destination of the trip. The estimated coefficients are $\hat{\beta} = 0.216$ for the Manhattan instance and $\hat{\beta} = 0.205$ for the four eastern boroughs. We define the driving durations of all arcs in the 2009 instances as the predictions of this model.

The 2018 \gls{TLC} data does not provide the latitude and longitude coordinates of the trips' origins and destinations. Instead, it provides the \gls{TLC} zone of each trip's origin and destination. We use this information to create a coarse-grained representation of the street network. We define the set of vertices as the set of all 263 \gls{TLC} zones in the four eastern boroughs. We define the set of arcs such that there is one pair of arcs between each pair of vertices whose zones border each other. Additionally, we manually add arcs to represent the bridges and tunnels that cross the East River and the Bronx River. We define the driving duration of each arc ($u,v$) as the median duration of all trips between zones $u$ and $v$.


\subsection{Demand distribution}
\label{sec:demand}

A demand realization at time $t$ is a customer request given by the vector \[b = (o_b, d_b, n_b, \taur_b, \taup_b, f_b)^{\intercal}.\]
The 2009 \gls{TLC} data contains the latitude and longitude coordinates of the origin $o_b$ and destination $d_b$ of every trip $b$. The 2018 \gls{TLC} data provides the \gls{TLC} zones of the origin and destination. For both years, the \gls{TLC} data contains the number of passengers $n_b$ and the trip fare $f_b$ of each request, as well as the pickup time, which we use as the time at which the \gls{rhs} receives the request. The latest pickup time $\taup_b = T$, i.e., all requests must be picked up before the end of the planning horizon. We assume that the latest response time $\taur_b$ is \qty{5}{\minute} after the request occurs.

We define a sample path as a set of demand realizations that occur between time $t=0$ and $t=T$. We associate each street network with a set of of sample paths that we derive from its respective \gls{TLC} data (2009 and 2018) as follows. Similarly to \citet{kullman:2022}, we use the data of March and April, a span of time that contains no public holidays. We take several steps to clean the data. We only consider the trips on Mondays to Fridays to reduce the variance. From these, we keep only those trips whose origin and destination are distinct, distance is between \qty{16.09}{\meter} and \qty{80}{\kilo\meter}, duration is greater than \qty{1}{min}, passenger count is positive, and whose fare is not below the base fare for taxis in \gls{NYC}. We then remove the trips with the 5\% largest fares, which tend to be outliers. Finally, we define the speed of a trip as its distance divided by its duration, and we remove all trips that are either slower than the average taxi speed in the most crowded part of \gls{NYC} or faster than the maximum speed limit on any street in \gls{NYC}, according to the \citet{nyc-mobility-report-2018}. 

\newcommand{\scalefactor}{$1\,\%$}
The \gls{TLC} data contains millions of requests. To keep the computational effort manageable, we define the number of requests in a sample path as \scalefactor{} of the trips in a randomly chosen Tuesday, Thursday or Friday in March or April of the year of the instance. We choose these trips randomly, as samples from a uniform distribution, separately for each sample path. To keep the ratio between requests and vehicles realistic, we define the fleet size as the number of simultaneously active taxis, scaled down by \scalefactor{} again. We gather taxi count estimates from \citet{toddwschneider}. All in all, our sampling method closely resembles that of \citet{kullman:2022}, except for the request count, which in their case is a parameter and constant across sample paths, but in our case is sampled from the \gls{TLC} data and can be different in every sample path. 

\subsection{Implementation details}
\label{sec:implementation}
We train the \gls{vfa} policy on 3000  sample paths, and we test both policies on 30 sample paths. For each instance, we create 3000 training sample paths and 30 test sample paths and reuse them across problem settings and policies. Similarly, we sample the initial state, which is the location and driving range of each vehicle, once per sample path from a uniform distribution. Every vehicle starts empty. As a result, all \gls{rhs} operators face exactly the same conditions in terms of initial state and demand, which enables a fair comparison.

We use $\theta = 0.1$ as the charging threshold of the \gls{pm} policy. This means that the \gls{pm} policy recharges a vehicle if and only if it has less than or equal to $10\,\%$ range left. We found that this threshold, which is the same as the one \citet{alkanj:2020} use, provides the best performance based on preliminary experiments. 

The \gls{vfa} policy we implement uses a look-up table with monotonic properties as described in \Cref{sec:monotonicity} of this online supplement to approximate the value of each vehicle attribute. 
We 
aggregate the driving range $l_a$ into nine equidistant levels, the capacity $n_a$ into two levels (vehicle is empty or not), and we aggregate the actionable time $t_a$ into 288 time levels, each of which corresponds to five minutes. This does not affect the decision epoch. We do not aggregate the location dimensions, i.e., the vehicle's origin $o_a$ and destination $d_a$. We use generalized harmonic learning rate decay with the parameter $a=300$ \citep{powell:2011}. For the demand attributes, instead of a lookup table, we use the estimate $\bar{v}_{tb}^D=0$ if $t^\text{res}_b \leq t$, and $\bar{v}_{tb}^D=0.9 f_b$ otherwise.


We implement our numerical experiments in Julia 1.10 using JuMP version 1.20 \citep{jump:2017} with the solver Gurobi 11.0 \citep{gurobi2024} to model and solve the LP and MIP problems. We performed the simulations on Intel Xeon 8468 Sapphire CPUs using up to \qty{450}{\giga\byte} of memory, depending on the instance. The code developed for this paper will be made available upon publication.

\FloatBarrier
\subsection{Tables}
\label{sec:tables}
The two tables in this section contain detailed results of our numerical experiments.
In Table \ref{tab:rfr_results}, we report the average \gls{rfr}, the median, the \gls{iqr}, and the \gls{me} for the 95\% confidence interval given by $1.96 \sigma\sqrt{N}^{-1}$ where $\sigma$ is the standard deviation of the cumulative rewards and $N$ is the number of sample paths in the test set.
In Table \ref{tab:abs_reward_results}, we report the same statistics regarding the absolute rewards.

\begin{table}[hbtp]
    \centering
    \resizebox{\textwidth}{!}{
        \begin{tabular}{@{}llllrrrrlrrrr@{}}
            
            \toprule \multirow{2}{*}{Instance} & \multirow{2}{*}{Pooling-enabled} & \multirow{2}{*}{Fleet type}
            &  & \multicolumn{4}{c}{PM Reward (\$)} & ~ & \multicolumn{4}{c}{VFA Reward (\$)} \\ \cmidrule{5-8} \cmidrule{10-13}
            &  &  &  & \multicolumn{1}{c}{Mean} & \multicolumn{1}{c}{Median} & \multicolumn{1}{c}{\acrshort{iqr}} & \multicolumn{1}{c}{\acrshort{me}} &  & \multicolumn{1}{c}{Mean} & \multicolumn{1}{c}{Median} & \multicolumn{1}{c}{\acrshort{iqr}} & \multicolumn{1}{c}{\acrshort{me}} \\
            \midrule
            \multirow[t]{6}{*}{1} & \multirow[t]{3}{*}{false} & \acrshort{EV} (\acrshort{DCFC}) &  & 18990 & 19017 & 1296 & 342 &  & 20508 & 20771 & 1456 & 365 \\
            &  & \acrshort{EV} (\acrshort{l2c}) &  & 13426 & 13465 & 836 & 221 &  & 18949 & 19065 & 1048 & 291 \\
            &  & \acrshort{ICE} &  & 19170 & 19459 & 1711 & 379 &  & 20408 & 20820 & 1808 & 392 \\
            
            & \multirow[t]{3}{*}{true} & \acrshort{EV} (\acrshort{DCFC}) &  & 19882 & 19998 & 1459 & 363 &  & 21240 & 21403 & 2022 & 534 \\
            &  & \acrshort{EV} (\acrshort{l2c}) &  & 12896 & 12884 & 1070 & 244 &  & 20741 & 20860 & 1907 & 506 \\
            &  & \acrshort{ICE} &  & 20272 & 20497 & 921 & 436 &  & 21146 & 21298 & 1991 & 534 \\
            \midrule 
            \multirow[t]{6}{*}{2} & \multirow[t]{3}{*}{false} & \acrshort{EV} (\acrshort{DCFC}) &  & 22382 & 22226 & 1410 & 395 &  & 24347 & 24006 & 2093 & 573 \\
            &  & \acrshort{EV} (\acrshort{l2c}) &  & 16589 & 16650 & 788 & 227 &  & 21812 & 21874 & 834 & 282 \\
            &  & \acrshort{ICE} &  & 22439 & 22267 & 1407 & 446 &  & 24326 & 24218 & 1570 & 494 \\
            
            & \multirow[t]{3}{*}{true} & \acrshort{EV} (\acrshort{DCFC}) &  & 23415 & 23352 & 1601 & 406 &  & 24305 & 23990 & 1882 & 541 \\
            &  & \acrshort{EV} (\acrshort{l2c}) &  & 15840 & 15787 & 889 & 214 &  & 23423 & 23077 & 1699 & 518 \\
            &  & \acrshort{ICE} &  & 23579 & 23581 & 1720 & 430 &  & 24263 & 23893 & 1685 & 540 \\
            \midrule 
            \multirow[t]{6}{*}{3} & \multirow[t]{3}{*}{false} & \acrshort{EV} (\acrshort{DCFC}) &  & 12349 & 12335 & 466 & 198 &  & 14849 & 14829 & 1006 & 330 \\
            &  & \acrshort{EV} (\acrshort{l2c}) &  & 8249 & 8286 & 688 & 162 &  & 13733 & 13754 & 644 & 198 \\
            &  & \acrshort{ICE} &  & 12469 & 12472 & 799 & 201 &  & 13700 & 13794 & 526 & 202 \\
            
            & \multirow[t]{3}{*}{true} & \acrshort{EV} (\acrshort{DCFC}) &  & 12617 & 12685 & 539 & 183 &  & 14895 & 14826 & 1252 & 370 \\
            &  & \acrshort{EV} (\acrshort{l2c}) &  & 7719 & 7641 & 698 & 171 &  & 13599 & 13736 & 531 & 208 \\
            &  & \acrshort{ICE} &  & 12973 & 13094 & 472 & 203 &  & 14938 & 14919 & 1261 & 387 \\
            \midrule 
            \multirow[t]{6}{*}{4} & \multirow[t]{3}{*}{false} & \acrshort{EV} (\acrshort{DCFC}) &  & 16700 & 16736 & 819 & 220 &  & 20191 & 20126 & 1525 & 289 \\
            &  & \acrshort{EV} (\acrshort{l2c}) &  & 11035 & 11006 & 487 & 136 &  & 18668 & 18698 & 696 & 191 \\
            &  & \acrshort{ICE} &  & 16872 & 16648 & 793 & 243 &  & 20288 & 20164 & 1063 & 242 \\
            
            & \multirow[t]{3}{*}{true} & \acrshort{EV} (\acrshort{DCFC}) &  & 17254 & 17227 & 536 & 152 &  & 20824 & 20538 & 1594 & 357 \\
            &  & \acrshort{EV} (\acrshort{l2c}) &  & 10394 & 10305 & 484 & 190 &  & 19576 & 19438 & 1005 & 225 \\
            &  & \acrshort{ICE} &  & 17636 & 17537 & 811 & 184 &  & 20990 & 20900 & 1791 & 375 \\
            \bottomrule
        \end{tabular}

    }
    \caption{Total reward statistics.}
    \label{tab:abs_reward_results}
\end{table}

\begin{table}[htbp]
    \centering
    \resizebox{\textwidth}{!}{%
        \begin{tabular}{@{}llllrrrrlrrrr@{}}
            \toprule \multirow{2}{*}{Instance} & \multirow{2}{*}{Pooling enabled} & \multirow{2}{*}{Fleet type}
            &  & \multicolumn{4}{c}{PM \acrshort{rfr} (\%)} & ~ & \multicolumn{4}{c}{VFA \acrshort{rfr} (\%)} \\ \cmidrule{5-8} \cmidrule{10-13}
            &  &  &  & \multicolumn{1}{c}{Mean} & \multicolumn{1}{c}{Median} & \multicolumn{1}{c}{\acrshort{iqr}} & \multicolumn{1}{c}{\acrshort{me}} &  & \multicolumn{1}{c}{Mean} & \multicolumn{1}{c}{Median} & \multicolumn{1}{c}{\acrshort{iqr}} & \multicolumn{1}{c}{\acrshort{me}} \\
            \midrule
            \multirow[t]{6}{*}{1} & \multirow[t]{3}{*}{false} & \acrshort{EV} (\acrshort{DCFC}) &  & 80.724 & 80.581 & 4.224 & 1.027 &  & 85.808 & 86.620 & 4.199 & 0.989 \\
            &  & \acrshort{EV} (\acrshort{l2c}) &  & 57.195 & 56.895 & 5.664 & 1.537 &  & 79.366 & 79.597 & 4.334 & 1.091 \\
            &  & \acrshort{ICE} &  & 81.372 & 81.458 & 3.490 & 1.009 &  & 85.249 & 84.947 & 3.419 & 0.825 \\
            
            & \multirow[t]{3}{*}{true} & \acrshort{EV} (\acrshort{DCFC}) &  & 85.611 & 85.672 & 2.093 & 0.823 &  & 89.331 & 89.204 & 1.154 & 0.311 \\
            &  & \acrshort{EV} (\acrshort{l2c}) &  & 55.438 & 55.245 & 4.870 & 1.369 &  & 87.419 & 87.186 & 1.768 & 0.395 \\
            &  & \acrshort{ICE} &  & 86.057 & 86.155 & 2.822 & 1.010 &  & 89.027 & 89.044 & 1.316 & 0.290 \\
            \midrule 
            \multirow[t]{6}{*}{2} & \multirow[t]{3}{*}{false} & \acrshort{EV} (\acrshort{DCFC}) &  & 83.358 & 83.016 & 3.515 & 0.767 &  & 89.280 & 89.280 & 0.460 & 0.154 \\
            &  & \acrshort{EV} (\acrshort{l2c}) &  & 61.867 & 62.399 & 5.535 & 1.331 &  & 79.831 & 79.764 & 4.066 & 1.117 \\
            &  & \acrshort{ICE} &  & 83.629 & 83.463 & 2.616 & 0.848 &  & 89.113 & 89.444 & 1.132 & 0.435 \\
            
            & \multirow[t]{3}{*}{true} & \acrshort{EV} (\acrshort{DCFC}) &  & 88.160 & 88.237 & 2.264 & 0.478 &  & 89.443 & 89.422 & 0.609 & 0.181 \\
            &  & \acrshort{EV} (\acrshort{l2c}) &  & 59.856 & 59.589 & 4.476 & 1.376 &  & 86.471 & 86.530 & 0.993 & 0.324 \\
            &  & \acrshort{ICE} &  & 88.660 & 89.087 & 2.070 & 0.421 &  & 89.080 & 89.074 & 0.639 & 0.175 \\
            \midrule 
            \multirow[t]{6}{*}{3} & \multirow[t]{3}{*}{false} & \acrshort{EV} (\acrshort{DCFC}) &  & 76.639 & 76.776 & 3.234 & 1.215 &  & 91.007 & 91.222 & 1.411 & 0.348 \\
            &  & \acrshort{EV} (\acrshort{l2c}) &  & 51.258 & 51.229 & 5.545 & 1.616 &  & 84.188 & 84.098 & 3.370 & 1.152 \\
            &  & \acrshort{ICE} &  & 77.411 & 76.903 & 3.843 & 1.309 &  & 83.781 & 84.408 & 3.794 & 1.051 \\
            
            & \multirow[t]{3}{*}{true} & \acrshort{EV} (\acrshort{DCFC}) &  & 79.034 & 78.778 & 3.820 & 1.077 &  & 91.596 & 91.607 & 0.732 & 0.209 \\
            &  & \acrshort{EV} (\acrshort{l2c}) &  & 48.251 & 47.339 & 5.619 & 1.452 &  & 83.736 & 83.883 & 3.615 & 0.997 \\
            &  & \acrshort{ICE} &  & 81.032 & 81.267 & 3.577 & 1.105 &  & 91.959 & 92.100 & 1.078 & 0.247 \\
            \midrule 
            \multirow[t]{6}{*}{4} & \multirow[t]{3}{*}{false} & \acrshort{EV} (\acrshort{DCFC}) &  & 73.798 & 72.832 & 3.535 & 0.929 &  & 88.318 & 88.402 & 1.567 & 0.454 \\
            &  & \acrshort{EV} (\acrshort{l2c}) &  & 48.723 & 49.351 & 2.872 & 0.921 &  & 81.839 & 81.998 & 4.296 & 1.017 \\
            &  & \acrshort{ICE} &  & 74.625 & 75.098 & 2.328 & 0.844 &  & 88.891 & 89.005 & 2.805 & 0.662 \\
            
            & \multirow[t]{3}{*}{true} & \acrshort{EV} (\acrshort{DCFC}) &  & 77.213 & 77.093 & 3.486 & 0.908 &  & 91.215 & 91.230 & 0.621 & 0.236 \\
            &  & \acrshort{EV} (\acrshort{l2c}) &  & 46.072 & 46.350 & 2.817 & 0.910 &  & 85.891 & 86.138 & 3.365 & 0.760 \\
            &  & \acrshort{ICE} &  & 78.876 & 78.779 & 2.979 & 0.658 &  & 91.879 & 91.963 & 0.788 & 0.241 \\
            \bottomrule
        \end{tabular}
    }
    \caption{\Acrfull{rfr} statistics.}
    \label{tab:rfr_results}
\end{table}

\FloatBarrier
\end{appendix}

\end{document}